\documentclass[a4paper,11pt,twoside]{article}
\RequirePackage[a4paper,head=1cm]{geometry}
\usepackage[utf8x]{inputenc}
\usepackage{textcomp}
\usepackage{amsmath,amsfonts,verbatim,afterpage,euscript,mathrsfs,amssymb,empheq}
\usepackage{amsthm}
\usepackage{dsfont}
\usepackage{hyperref}
\usepackage{pst-fill,pst-grad}
\usepackage{textcomp}
\usepackage{multicol}
\usepackage{amsmath}
\usepackage{amssymb}
\usepackage{hyperref}
\usepackage{dsfont}

\newtheorem{Theoreme}{Theorem}
\newtheorem{Lemme}{Lemma}[section]

\newtheorem{Proposition}{Proposition}[section]
\newtheorem{Remarque}{\bf Remark}
\newcommand{\mysection}{\setcounter{equation}{0} \section}

\def \vu{\vec{u}}
\def \vpsi{\vec{\psi}}

\def \vn{\vec{\nabla}}

\title{\bf Global mild solutions for a transport-diffusion equation with a rough drift.} 
\author{Diego Chamorro\footnote{Laboratoire de Math\'ematiques et Mod\'elisation d'Evry (LaMME) - UMR 8071. Universit\'e d'Evry Val d'Essonne, 23 Boulevard de France, 91037 Evry Cedex, France. email: \textit{diego.chamorro@univ-evry.fr}}, Anca-Nicoleta Marcoci\footnote{Department of Mathematics and Computer Science. Technical University of Civil Engineering, Bucharest, Bld. Lacul Tei, no. 124, sector 2. Romania. email: \textit{anca.marcoci@utcb.ro}}, Liviu-Gabriel Marcoci\footnote{Department of Mathematics and Computer Science. Technical University of Civil Engineering, Bucharest, Bld. Lacul Tei, no. 124, sector 2. Romania.  email: \textit{liviu.marcoci@utcb.ro}}.} 
\begin{document} 
%%%%%%%%%%%%%%%%%%%%%%%%%%%%%%%%%%%%%%%%%%%%%%%%%%%
\maketitle 
%%%%%%%%%%%%%%%%%%%%%%%%%%%%%%%%%%%%%%%%%%%%%%%%%%%
%%%%%%%%%%%%%%%%%%%%%%%%%%%%%%%%%%%%%%%%%%%%%%%%%%%
\begin{scriptsize}
\abstract{ \noindent We construct here global mild solutions in a critical setting for a class of transport-diffusion equations with a drift term that involves rough Calder\'on-Zygmund operators.}\\

\noindent \textbf{MSC (2020) Primary: 35A01; Secondary: 35B33}
\end{scriptsize}
%%%%%%%%%%%%%%%%%%%%%%%%%%%%%%%%%%%%%%%%%%%%%%%%%%%
%\tableofcontents 

%%%%%%%%%%%%%%%%%%%%%%%%%%%%%%%%%%%%%%%%%%%%%%%%%%%
%%%%%%%%%%%%%%%%%%%%%%%%%%%%%%%%%%%%%%%%%%%%%%%%%%%
%%%%%%%%%%%%%%%%%%%%%%%%%%%%%%%%%%%%%%%%%%%%%%%%%%%
\mysection{Introduction}
In this article we will study some variations of a nonlinear transport-diffusion equation in $\mathbb{R}^n$ with $n\geq 2$ where the drift term is given by a vector of rough singular integral operators. We are particularly interested in the existence of global in time mild solutions in a critical resolution space associated to some (small) initial data. In this sense, the most representative equation that will be studied here is the following:
\begin{equation}\label{EquationIntro}
\begin{cases}
\partial_t \theta(t,x) = \Delta \theta(t,x) + (\mathbb{T}_{[\theta]}\cdot \nabla \theta)(t,x)+f(t,x), \qquad (t,x)\in [0,+\infty[\times \mathbb{R}^n,\\[3mm]
\theta(0,x)=\theta_0(x), \qquad x\in \mathbb{R}^n,
\end{cases}
\end{equation}
where $\theta:[0,+\infty[\times \mathbb{R}^n\longrightarrow \mathbb{R}$ is the unknown, $f:[0,+\infty[\times \mathbb{R}^n\longrightarrow \mathbb{R}$ is a given external force, $\theta_0:\mathbb{R}^n\longrightarrow \mathbb{R}$ is an initial data and $\mathbb{T}_{[\theta]}=(T_1(\theta), \cdots, T_n(\theta))$ is a vector composed of $n$ rough singular integral operators $(T_k)_{1\leq k\leq n}$ which are all defined by the expression 
\begin{equation}\label{Def_Operator}
T_k(\phi)(x)=p.v.\int_{\mathbb{R}^n}\frac{\Omega_k(y/ |y|)}{|y|^n}\phi(x-y)dy,\qquad 1\leq k\leq n,
\end{equation}
for any locally integrable function $\phi:\mathbb{R}^n\longrightarrow \mathbb{R}$, and where the kernels $\Omega_k:\mathbb{S}^{n-1}\longrightarrow \mathbb{R}$ are such that 
$\Omega_k \in L^1(\mathbb{S}^{n-1})$, $\displaystyle{\int_{\mathbb{S}^{n-1}}\Omega_k \ d\sigma=0}$ and $\Omega_k \in L^{\rho}(\mathbb{S}^{n-1})$ with $1<\rho<n$. Note that these operators fall outside the ``classical'' setting of singular integral operators of convolution type (see Section 4.4 of the book \cite{Grafakos}) and in this sense we will say that the vector $\mathbb{T}_{[\theta]}$ represents a \emph{rough} drift. Some properties of this general family of rough singular operators were studied in \cite{Hoang1}, \cite{Hoang} and \cite{Li} (although with different conditions over the kernel function $\Omega$). Remark also that since the drift depends on the function $\theta$, we are thus dealing with a nonlinear transport-diffusion equation. Note finally that the operators $T_k$ act only in the space variable and no interaction is asked in the time variable.\\

The system (\ref{EquationIntro}) can be seen as a generalization of a class of partial differential equations which came from fluid dynamics. For example, if $n=2$ and if we set $\mathbb{T}_{[\theta]}=(-R_2(\theta), R_1(\theta))$, where $R_j$ with $j=1,2$ are the classical Riesz transforms (which are singular integrals of convolution type), we obtain the 2D quasi-geostrophic equation which was intensively studied in a fractional setting, \emph{i.e.} when the Laplace operator $\Delta$ is replaced by its fractional power $-(-\Delta)^{\frac{\alpha}{2}}$ for $0<\alpha<2$:
\begin{equation}\label{SQG}
\partial_t \theta = -(-\Delta)^{\frac{\alpha}{2}}\theta + (\mathbb{T}_{[\theta]}\cdot \nabla \theta).
\end{equation}
See \emph{e.g.} \cite{Caffarelli}, \cite{ChamorroMenozzi}, \cite{Cordoba} and the references therein for more details about this equation. \\

We can also consider the 3D Navier-Stokes equation $\partial_t\vu=\Delta\vu-(\vu\cdot \vn)\vu-\vn p$,  with $div(\vu)=0$, which can be rewritten, after the application of the Leray projector $\mathbb{P}(\vpsi)=\vpsi+\vn (-\Delta)^{-1}\vn\cdot \vpsi$, as 
\begin{equation}\label{NS}
\partial_t\vu=\Delta\vu-\mathbb{P}((\vu\cdot \vn)\vu).
\end{equation}
Note that this is also a nonlinear transport-diffusion equation that involves singular integral operators since the Riesz transforms are intrinsically used in the definition of the Leray projector (although the form of this equation is slightly different from the one of the equation (\ref{EquationIntro}) since the singular integral operator $\mathbb{P}$ acts on the whole nonlinear term $(\vu \cdot \vn)\vu$). See the book \cite{Chamorro_Livre} for essential facts about the Navier-Stokes equations and see the book \cite{PGLR} for a more advanced study of this system. The equations (\ref{SQG}) and (\ref{NS}) contain many open problems that can be studied in very different directions (existence of mild or weak solutions, energy inequalities, regularity properties, uniqueness, etc.).\\

The starting point to construct \emph{mild} solutions for the system (\ref{EquationIntro}) is given by the integral representation formula (also known as the Duhamel formula):
\begin{equation}\label{EquationIntro_Integrale}
\theta(t,x)=\mathfrak{g}_t\ast \theta_0(x)+\int_{0}^{+\infty}\mathfrak{g}_{t-s}\ast (\mathbb{T}_{[\theta]}\cdot \nabla \theta)(s,x)ds+\int_{0}^{+\infty}\mathfrak{g}_{t-s}\ast f(s,x)ds,
\end{equation}
and with this integral representation formula we can easily apply the classical Banach-Picard fixed point algorithm (see \cite[Théor\`eme 4.1.1]{Chamorro_Livre} for a proof) to obtain a mild solution to the previous problem: indeed, consider $\theta_0:\mathbb{R}^n\longrightarrow \mathbb{R}$ an initial data that belongs to a Banach space $E_0(\mathbb{R}^n)$ and an external force $f:[0,+\infty[\times\mathbb{R}^n\longrightarrow \mathbb{R}$ that belongs to a Banach space $E_f([0,+\infty[\times\mathbb{R}^n)$. If
$E([0,+\infty[\times\mathbb{R}^n)$ is a Banach space (usually called the \emph{resolution space}) such that we have the controls 
\begin{eqnarray}
&&\|\mathfrak{g}_t\ast \theta_0\|_{E}\leq C_0\|\theta_0\|_{E_0},\qquad \qquad \left\|\int_{0}^{+\infty}\mathfrak{g}_{t-s}\ast fds\right\|_{E}\leq C_f\|f\|_{E_f},\label{Estimations_Donnees}\\[3mm]
&&\left\|\int_{0}^{+\infty}\mathfrak{g}_{t-s}\ast (\mathbb{T}_{[\theta]}\cdot \nabla \theta)ds\right\|_{E}\leq C_B\|\theta\|_E\|\theta\|_E,\label{Estimation_Nonlineaire}
\end{eqnarray}
and if we have the relationship $C_0\|\theta_0\|_{E_0}+C_f\|f\|_{E_f}<\frac{1}{4C_B}$, then by applying the Banach-Picard contraction principle we obtain a solution $\theta\in E$ that is a \emph{mild} solution of the equation (\ref{EquationIntro_Integrale}).\\

As we can see, we can construct mild solutions to the system (\ref{EquationIntro_Integrale}) as long as we have good functional inequalities that can lead us to estimates of the form (\ref{Estimations_Donnees}) and (\ref{Estimation_Nonlineaire}). In this sense, although very similar to the system (\ref{EquationIntro}), the equations (\ref{SQG}) and (\ref{NS}) are more ``flexible'', since the Riesz transforms $(R_j)_{1\leq j\leq n}$ are bounded in \emph{all} the usual spaces (for example in the Lebesgue spaces $L^p(\mathbb{R}^n)$ with $1<p<+\infty$, we have $\|R_j(\theta)\|_{L^p}\leq C\|\theta\|_{L^p}$). However, in the case of the rough singular integral operators $(T_k)_{1\leq k\leq n}$ considered in the expression (\ref{Def_Operator}) above and associated to a kernel $\Omega_k\in L^{\rho}(\mathbb{S}^{n-1})$ with $1<\rho<n$, we do not have, to the best of our knowledge, such flexibility and this will introduce some rigidity in the functional framework that will be used to obtain \emph{mild} solutions. Indeed, we only know how to prove the generic estimate $\|T_k(\theta)\|_{L^{\mathfrak{s}}}\leq C\|\nabla \theta\|_{L^q}$ with $\mathfrak{s}=\frac{nq}{n-q}$ and $1<\frac{\rho n}{\rho n+\rho-n}\leq q<n$ (see  the Lemma \ref{LemmaRough} below) and this implies in particular that the ``price to pay'' to obtain a control of the quantity $\|T_k(\theta)\|_{L^{\mathfrak{s}}}$ is to consider a $L^q$-norm of the gradient $\nabla \theta$, but not of the function $\theta$, and we don't know if a better estimate in the range of Lebesgue spaces is available for this type of rough operators that does not involves a gradient. This lack of flexibility in the functional estimates will be a determining factor when choosing functional spaces in order to close the fixed point argument.\\

As mentioned above, we are interested here in \emph{global in time} mild solutions of the problem (\ref{EquationIntro_Integrale}) and one particular way to achieve this is to consider a functional framework that is adapted to the dilation structure of the equation. Indeed, if $\theta(t,x)$ is a solution of the equation (\ref{EquationIntro_Integrale}) associated to an initial data $\theta_0(x)$ and to an external force $f(t,x)$, and if we define $\theta_\lambda(t,x)=\lambda\theta(\lambda^2 t, \lambda x)$,  for some real parameter $\lambda>0$, then it is not difficult to see that the function $\theta_\lambda(t,x)$ is also a solution of the problem (\ref{EquationIntro_Integrale}) with an initial data $\theta_{\lambda 0}(t,x)=\lambda \theta_0(\lambda x)$ and with an external force $f_{\lambda}(t,x)=\lambda^3 f(\lambda^2 t,\lambda x)$ (note in particular that this scaling property with respect to the time and space dilation is the same as the one for the Navier-Stokes system (\ref{NS})). Thus, if the functional spaces $E_0$, $E_f$ and $E$ satisfy the scaling invariances 
\begin{equation}\label{Invariance}
\|\theta_{\lambda 0}\|_{E_0}=\|\theta\|_{E_0}, \qquad \|f_{\lambda}\|_{E_f}=\|f\|_{E_f},\quad \mbox{and}\quad \|\theta_\lambda\|_E=\|\theta\|_{E},
\end{equation}
the constants $C_0, C_f$ and $C_B$ in (\ref{Estimations_Donnees}) and (\ref{Estimation_Nonlineaire}) may not depend on the time variable: we will thus obtain global in time mild solutions under an inevitable smallness assumption.\\

We will thus study in this article global mild solutions in a critical framework for the system (\ref{EquationIntro}) as well as some variations of this nonlinear transport-diffusion equation.\\

In our first result, we consider an example of scaling invariant functional spaces for the initial data $E_0$, the external force $E_f$ and for the resolution space $E$ that will lead us to global in time mild solutions for the problem (\ref{EquationIntro}):
%%%%%%%%%%%%%%%%%%%%%%%%%%%%%%%%%%%%%%%%%%%%%%%%%%%
\begin{Theoreme}[\bf Global Mild Solutions for the system (\ref{EquationIntro})]\label{Theoreme1}
Consider a drift vector $\mathbb{T}_{[\theta]}$ conformed by rough singular operators $T_k$ of the form given in (\ref{Def_Operator}) associated to kernels $\Omega_k \in L^\rho(\mathbb{S}^{n-1})$ with $1<\rho<n$.\\

\noindent Let $\theta_0:\mathbb{R}^n\longrightarrow \mathbb{R}$ be an initial data that belongs to the homogeneous Besov space $E_0(\mathbb{R}^n)=\dot{B}^{\frac{n-q}{q}, q}_{\infty}(\mathbb{R}^n)$: 
\begin{equation}\label{Espace_ConditionInitiale}
\|\theta_0\|_{E_0}=\underset{t>0}{\sup}\; t^{\frac{2q-n}{2q}}\|\mathfrak{g}_{t}\ast \theta_0\|_{\dot{W}^{1,q}}=\|\theta_0\|_{\dot{B}^{\frac{n-q}{q}, q}_{\infty}}<+\infty, 
\end{equation}
with $1<\frac{\rho n}{\rho n+\rho-n}<q<n<2q$ and let $f:[0,+\infty[\times \mathbb{R}^n\longrightarrow \mathbb{R}$ be an external force that satisfies
\begin{equation}\label{Espace_ForceExterieure}
\|f\|_{E_f}= \underset{t>0}{\sup}\;t^{\frac{3}{2}-\frac{n}{2\varrho}}\|f(t,\cdot)\|_{L^\varrho}<+\infty, \quad \mbox{with } 1<\varrho<n<3\varrho \mbox{ and } \frac{1}{q}-\frac{1}{n}<\frac{1}{\varrho}<\frac{1}{q}. 
\end{equation}
If the quantity $\|\theta_0\|_{E_0}+\|f\|_{E_f}$ is small enough, then there exists a global in time mild solution $\theta(\cdot, \cdot)$ of the equation (\ref{EquationIntro}) such that 
\begin{equation}\label{Espace_Resolution}
\|\theta\|_{E}=\underset{t>0}{\sup}\; t^{\frac{2q-n}{2q}}\| \theta(t, \cdot)\|_{\dot{W}^{1,q}}<+\infty.
\end{equation}
\end{Theoreme}
%%%%%%%%%%%%%%%%%%%%%%%%%%%%%%%%%%%%%%%%%%%%%%%%%%%
\noindent Some remarks are in order here. Note that the functional spaces $E_0(\mathbb{R}^n)=\dot{B}^{\frac{n-q}{q}, q}_{\infty}(\mathbb{R}^n)$, $E_f([0,+\infty[\times \mathbb{R}^n)$ and the resolution space $E([0,+\infty[\times \mathbb{R}^n)$ are indeed scaling invariant with respect of the dilation structure of the equation -in the sense that we have the identities (\ref{Invariance})- and we will show in Section \ref{Secc_ProofTh1} that the constants $C_0, C_f, C_B$ that appear in the estimates (\ref{Estimations_Donnees}) and (\ref{Estimation_Nonlineaire}) do not depend on the time variable: we will thus obtain, in a very natural manner, global in time mild solutions for the system (\ref{EquationIntro_Integrale}). Note also that the lower condition $1<\frac{\rho n}{\rho n+\rho-n}<q$ is related to the information available over the kernels $\Omega_k\in L^\rho(\mathbb{S}^{n-1})$ (with $1<\rho<n$) while the condition $q<n<2q$ is essentially technical. Thus, the values of the parameter $q$ that govern the initial data as well as the resolution space and the external force are driven by the information available over the kernels $\Omega_k$: if $\rho\to n$, these kernels are in some sense less rough and we have $\frac{\rho n}{\rho n+\rho-n}\to 1$ which gives us a wide interval of possibilities  for the parameter $q$, however if $\rho\to 1$ the operator is rougher and this forces $\frac{\rho n}{\rho n+\rho-n}\to n$, restricting the choice of the parameter $q$. As we can see, the values of the parameter $q$ are very sensitive to the information available over the kernels $\Omega_k$ and for this reason, in order to keep things simple, we asked the same information for all the kernels $\Omega_k$.\\

Remark now that we can also consider the more common (but smaller\footnote{We have the inequality $\|\theta_0\|_{\dot{B}^{\frac{n-q}{q}, q}_{\infty}}\leq \|\theta_0\|_{\dot{W}^{\frac{n-q}{q}, q}}$, from which we deduce the inclusion $\dot{W}^{\frac{n-q}{q}, q}(\mathbb{R}^n)\subset\dot{B}^{\frac{n-q}{q}, q}_{\infty}(\mathbb{R}^n)$.}) homogeneous Sobolev space $\dot{W}^{\frac{n-q}{q},q}(\mathbb{R}^n)$ for the initial data $\theta_0$ instead of the Besov space $\dot{B}^{\frac{n-q}{q}, q}_{\infty}(\mathbb{R}^n)$ without any change on the conclusion of the previous theorem. Note in particular that when $n=3$ and $q=2$ (and thus we should have $\Omega_k\in L^\rho(\mathbb{S}^{n-1})$ with $\frac{6}{5}<\rho<3$), we obtain for the initial data the Sobolev space $\dot{H}^{\frac{1}{2}}(\mathbb{R}^3)$ which was considered by Fujita \& Kato in \cite{FujitaKato} in the context of the Navier-Stokes equations.\\ 

However, if we are interested in considering the \emph{largest} functional space for the initial data $\theta_0$ then, due to the maximality of the homogeneous Besov spaces (see \cite{Meyer}), we should consider the Besov space $E_0(\mathbb{R}^n)=\dot{B}^{-1,\infty}_{\infty}(\mathbb{R}^n)$, but based on the work of Bourgain \& Pavlovic \cite{besov_ill_posed} in the Navier-Stokes equations, this functional space seems to be completely out of reach for the system (\ref{EquationIntro}) since the nonlinear term $(\mathbb{T}_{[\theta]}\cdot \nabla \theta)$ is far more rigid (due to the presence of the rough drift) than the term $(\vu\cdot \vn )\vu$. The study of a more general initial data $\theta_0$ that the one considered here will probably deserve a separated study and this constitutes a new open problem in the setting of the rough drift nonlinear equation (\ref{EquationIntro}).\\

In the previous theorem, we asked for the initial data the condition $\theta_0\in \dot{B}^{\frac{n-q}{q}, q}_{\infty}(\mathbb{R}^n)$ with $1<\frac{\rho n}{\rho n+\rho-n}<q<n<2q$. We will see now that, if we modify the equation (\ref{EquationIntro}), then it will be possible to consider an initial data in a negative regularity homogeneous space. Thus, our first attempt to modify the equation (\ref{EquationIntro}) so that we could consider an initial datum in a maximal homogeneous Besov space is the following:
\begin{equation}\label{EquationIntroDiv}
\begin{cases}
\partial_t \theta(t,x) =\Delta \theta(t,x) + div(\mathbb{T}_{[(-\Delta)^{-\frac{1}{2}}\theta]}\theta)(t,x)+f(t,x), \qquad (t,x)\in [0,+\infty[\times \mathbb{R}^n,\\[3mm]
\theta(0,x)=\theta_0(x), \qquad x\in \mathbb{R}^n.
\end{cases}
\end{equation}
If we compare the previous equation to the system (\ref{EquationIntro}), we can see two major modifications. The first one is related to the nonlinear term which is written in a divergence form: we have now $div(\mathbb{T}\theta)$ instead of $\mathbb{T}\cdot \nabla \theta$. Note that in fluid dynamics, a divergence-free property is usually asked, and therefore this modification is generally harmless. The second modification, stronger than the previous one, is a regularization of the rough singular integral drift $\mathbb{T}$ which is meant to compensate the lack of boundedness in Lebesgue spaces of the operators $T_k$: this will allow us to consider a more general initial data, but it will also change the scaling of the equation. Indeed, note now that if $\theta(t,x)$ is a solution of the equation (\ref{EquationIntroDiv}), associated to an initial data $\theta_0(x)$ and an external force $f(t,x)$, then for any $\nu>0$, the function $\theta_\nu(t,x)=\nu^2\theta(\nu^2t,\nu x)$ is also a solution of the equation  (\ref{EquationIntroDiv}), associated to an initial data $\theta_{\nu 0}(x)=\nu^2\theta_0(\nu x)$ and an external force $f_\nu(t,x)=\nu^4f(\nu^2t,\nu x)$. These relationships are related to the scaling invariance of the equation and they give a hint where to find critical functional spaces.\\

Our next result is the following one:
%%%%%%%%%%%%%%%%%%%%%%%%%%%%%%%%%%%%%%%%%%%%%%%%%%%
\begin{Theoreme}[\bf Global Mild Solutions for the system (\ref{EquationIntroDiv})]\label{Theoreme2}
Consider a drift vector $\mathbb{T}$ conformed by rough singular operators $T_k$ of the form given in (\ref{Def_Operator}) associated to kernels $\Omega_k \in L^\rho(\mathbb{S}^{n-1})$ with $1<\rho<n$.\\

\noindent Consider $\theta_0:\mathbb{R}^n\longrightarrow \mathbb{R}$ an initial data such that $\theta_0\in \dot{B}^{-\frac{2q-n}{q},q}_\infty(\mathbb{R}^n)$, with $1<\frac{\rho n}{\rho n+\rho-n}<q<n<2q$, \emph{i.e.}:
\begin{equation}\label{Espace_ConditionInitialeDiv}
\|\theta_0\|_{\dot{B}^{-\frac{2q-n}{q},q}_\infty}=\underset{t>0}{\sup}\; t^{\frac{2q-n}{2q}}\|\mathfrak{g}_{t}\ast \theta_0\|_{L^{q}}<+\infty,
\end{equation}
and let $f:[0,+\infty[\times \mathbb{R}^n\longrightarrow \mathbb{R}$ be an external force such that 
\begin{equation}\label{Espace_ForceExterieureDiv}
\|f\|_{\mathbb{E}_f}= \underset{t>0}{\sup}\;t^{(\frac{3}{2}-\frac{n}{2\varrho})} \|f(t,\cdot)\|_{\dot{W}^{-1,\varrho}}<+\infty, \end{equation}
with $1<\varrho<n<3\varrho$ and $\frac{1}{q}-\frac{1}{n}<\frac{1}{\varrho}<\frac{1}{q}$.\\ 

\noindent If the quantity $\|\theta_0\|_{ \dot{B}^{-\frac{2q-n}{q},q}_\infty}+\|f\|_{\mathbb{E}_f}$ is small enough, then  there exists a global in time mild solution $\theta$ of the equation (\ref{EquationIntroDiv}) such that 
\begin{equation}\label{Espace_ResolutionDiv}
\|\theta\|_{\mathbb{E}}=\underset{t>0}{\sup}\; t^{\frac{2q-n}{2q}}\| \theta(t, \cdot)\|_{L^{q}}<+\infty.
\end{equation}
\end{Theoreme}
%%%%%%%%%%%%%%%%%%%%%%%%%%%%%%%%%%%%%%%%%%%%%%%%%%%
\noindent We first remark that the functional spaces used in the previous result are indeed critical with respect to the scaling of the equation (\ref{EquationIntroDiv}) as we have the indentities
$$\|\theta_{\nu 0}\|_{\dot{B}^{-\frac{2q-n}{q},q}_\infty}=\|\theta_0\|_{\dot{B}^{-\frac{2q-n}{q},q}_\infty}, \quad \|f_{\nu}\|_{\mathbb{E}_f}=\|f\|_{\mathbb{E}_f},\quad\mbox{and}\quad \|\theta_{\nu}\|_{\mathbb{E}}=\|\theta\|_{\mathbb{E}},$$
and we will see that the constants $C_0, C_f$ and $C_B$ in the estimates (\ref{Estimations_Donnees}) and (\ref{Estimation_Nonlineaire}) will not depend on the time variable which will guarantee the existence of global in time mild solutions for the system (\ref{EquationIntroDiv}). Note now that, although we can now consider the Besov space $\dot{B}^{-\frac{2q-n}{q},q}_\infty(\mathbb{R}^n)$ for the initial data, due to the modification of the drift term made in the equation (\ref{EquationIntroDiv}) which introduces a different scaling, this Besov space $\dot{B}^{-\frac{2q-n}{q},q}_\infty(\mathbb{R}^n)$ can not be compared to the Besov space used in the Theorem \ref{Theoreme1}. \\

The modification introduced in the equation (\ref{EquationIntroDiv}) is quite strong but as expected it allows us to consider the negative regularity homogeneous Besov space $\dot{B}^{-\frac{2q-n}{q},q}_\infty(\mathbb{R}^n)$. However, the maximal homogeneous Besov space is in this case $\dot{B}^{-2,\infty}_\infty(\mathbb{R}^n)$ (recall that we have the embedding $\dot{B}^{-\frac{2q-n}{q},q}_\infty(\mathbb{R}^n) \subset \dot{B}^{-2,\infty}_\infty(\mathbb{R}^n)$) and despite of the regularization of the drift term this maximal functional space seems to be out of reach.\\

In our last result, we will modify the nonlinear drift term $\mathbb{T}_{[\theta]}$ in its inner structure  (but maintaining the same functions $\Omega_k$) and, for a parameter $0<\alpha<n$, we define now the operator $T_k^\alpha(\phi)$ by the expression
\begin{equation}\label{Def_OperatorAlpha}
T_k^\alpha(\phi)(x)=p.v.\int_{\mathbb{R}^n}\frac{\Omega_k(y/ |y|)}{|y|^{n-\alpha}}\phi(x-y)dy.
\end{equation}
Note that this is an operator of degree $-\alpha$ as we have $T_k^{\alpha}(\phi_\gamma)(x)=\gamma^{-\alpha}T_k^{\alpha}(\phi)(\gamma x)$ if $\phi_\gamma(x)=\phi(\gamma x)$ and thus some regularizing effect is to be expected: indeed, putting aside the function $\Omega_k$, the operator $T_k^\alpha$ can be seen as a generalization of the Riesz potential $I_\alpha$ which is given by the expression $I_\alpha(\phi)(x)=\displaystyle{c_\alpha\int_{\mathbb{R}^n}\frac{\phi(x-y)}{|y|^{n-\alpha}}dy}$ which has a regularization effect.\\

Now, if $\mathbb{T}^\alpha_{[\theta]}=(T_1^\alpha(\theta), \cdots, T_n^\alpha(\theta))$ is a vector of singular integral operators $T_k^\alpha$ of the form (\ref{Def_OperatorAlpha}), we will study the equation
\begin{equation}\label{EquationIntroAlpha}
\begin{cases}
\partial_t \theta(t,x) = \Delta \theta(t,x) + (\mathbb{T}_{[\theta]}^{\alpha}\cdot \nabla \theta)(t,x)+f(t,x), \qquad (t,x)\in [0,+\infty[\times \mathbb{R}^n,\\[3mm]
\theta(0,x)=\theta_0(x), \qquad x\in \mathbb{R}^n,
\end{cases}
\end{equation}
and in the Theorem \ref{Theoreme3} below we will obtain global in time mild solutions for this equation in an adapted critical framework. In particular, remark that if $\theta(t,x)$ is a solution of the problem (\ref{EquationIntroAlpha}) associated to an initial data $\theta_0(x)$ and an external force $f(t,x)$, then for all $\gamma>0$ the function $\theta_\gamma(t,x)=\gamma^{1+\alpha}\theta(\gamma^2t, \gamma x)$ is also a solution of the same problem associated to the initial data $\theta_{\gamma 0}(x)=\gamma^{1+\alpha}\theta_0(\gamma x)$ and to the external force $f_\gamma(t,x)=\gamma^{3+\alpha}f(\gamma^2t, \gamma x)$.\\

Following this critical framework we have the next existence result of global mild solutions:
%%%%%%%%%%%%%%%%%%%%%%%%%%%%%%%%%%%%%%%%%%%%%%%%%%%
\begin{Theoreme}[\bf Global Mild Solutions for the system (\ref{EquationIntroAlpha})]\label{Theoreme3}
For the equation (\ref{EquationIntroAlpha}), consider a drift vector $\mathbb{T}^\alpha_{\theta}$ with $0<\alpha<n-1$ conformed by rough singular operators $T^\alpha_k$ of the form given in (\ref{Def_OperatorAlpha}) associated to kernels $\Omega_k \in L^\rho(\mathbb{S}^{n-1})$ with $1<\rho<n$.\\

\noindent Let $\theta_0$ be an initial data such that 
\begin{equation}\label{Espace_ConditionInitialeAlpha}
\|\theta_0\|_{\mathcal{E}_0}=\underset{t>0}{\sup}\; t^{\frac{(2+\alpha)q-n}{2q}}\|\mathfrak{g}_{t}\ast \theta_0\|_{\dot{W}^{1,q}}=\|\theta_0\|_{\dot{B}^{\frac{n-(1+\alpha)q}{q},q}_\infty}<+\infty,
\end{equation}
with $1<\max\{\frac{\rho n}{\rho n+\rho-n}, \frac{n}{2+\alpha}\}<q<\frac{n}{1+\alpha}$ and let $f$ be an external force such that 
\begin{equation}\label{Espace_ForceExterieureAlpha}
\|f\|_{\mathcal{E}_f}=\underset{t>0}{\sup}\; t^{\frac{(3+\alpha)}{2}-\frac{n}{2\varrho}}\|f(t, \cdot)\|_{L^\varrho}<+\infty,
\end{equation}
with $\max\{1, \frac{n}{3+\alpha}\}<\varrho<\frac{n}{1+\alpha}$ and $\frac{1}{q}-\frac{1}{n}<\frac{1}{\varrho}<\frac{1}{q}$.\\

\noindent If the quantity $\|\theta_0\|_{\mathcal{E}_0}+\|f\|_{\mathcal{E}_f}$ is small enough, then there exists a global in time mild solution $\theta$ of the equation (\ref{EquationIntroAlpha}) such that 
\begin{equation}\label{Espace_ResolutionAlpha}
\|\theta\|_{\mathcal{E}}=\underset{t>0}{\sup}\; t^{\frac{(2+\alpha)q-n}{2q}}\| \theta(t, \cdot)\|_{\dot{W}^{1,q}}<+\infty.
\end{equation}
\end{Theoreme}
%%%%%%%%%%%%%%%%%%%%%%%%%%%%%%%%%%%%%%%%%%%%%%%%%%%
\noindent We first remark now that the functional framework used in this theorem is indeed adapted to the scaling of the equation (\ref{EquationIntroAlpha}) as we have the identities $\|\theta_{\gamma 0}\|_{\mathcal{E}_0}=\|\theta_{0}\|_{\mathcal{E}_0}$, $\|f_{\gamma}\|_{\mathcal{E}_f}=\|f\|_{\mathcal{E}_f}$ and $\|\theta_\gamma\|_{\mathcal{E}}=\|\theta_\gamma\|_{\mathcal{E}}$ and this allows us to obtain global in time solutions. Next note that if we set $\alpha\to 0$, then from Theorem \ref{Theoreme3} above we recover the results of the Theorem \ref{Theoreme1}, and this indicates a form of stability of our method to obtain mild solutions with respect to the ``regularization'' properties of the operators $T^\alpha_k$ defined in (\ref{Def_OperatorAlpha}). Let us mention that the key result that allows us to obtain suitable functional inequalities needed to close the fixed point argument is the following estimate 
$$\|T_k^\alpha(\phi)\|_{L^{\mathfrak{s}}}\leq C\|\nabla \phi\|_{L^q},$$
where $\mathfrak{s}=\frac{qn}{n-q(1+\alpha)}>1$ and $1<\frac{\rho n}{\rho n+\rho-n}<q<\frac{n}{1+\alpha}$ (see Proposition \ref{Proposition1} below). Note that the regularization effect of this estimate is only related to the choice of the parameters that fix the Lebesgue space in the left-hand side of the inequality above as we still have a gradient in the right-hand side. The previous estimate, which is proven in the Appendix \ref{Secc_ProofProposition1}, seems to be new in the context of rough singular integral operators. \\[2mm]

To conclude this section, let us make some few general remarks. First note that we do not claim any optimality for the parameters used in our previous theorems and perhaps a different study may lead to a wider range of indexes. Note now that the resolution spaces $E$, $\mathbb{E}$ or $\mathcal{E}$ introduced previously share a common structure: a Lebesgue or Sobolev $L^q$ norm in the space variable and a weighted $L^\infty$ norm in the time variable. Other resolution spaces may also be considered and we do not claim to be exhaustive. Finally, as it was mentioned above, the study of the largest critical functional space for an initial data (in order to obtain mild solutions) is sometimes a very challenging open problem and although there is a \emph{natural} functional framework given by the negative regularity homogeneous Besov spaces $\dot{B}^{-s,\infty}_\infty$, these critical spaces are not always accessible and this seems to be the case for the equations considered here where the main problem relies on suitable functional inequalities -even for the equation (\ref{EquationIntroDiv}) where the nonlinear drift was regularized. Indeed, to the best of our knowledge, for controlling the rough singular integral operators $T_k$ or $T_k^\alpha$, we only dispose (in the setting of Lebesgue spaces) of the estimates given in the Lemma \ref{LemmaRough} and in the Proposition \ref{Proposition1} which seems not enough to consider larger functional spaces for the initial data.\\

%%%%%%%%%%%%%%%%%%%%%%%%%%%%%%%%%%%%%%%%%%%%%%%%%%%
The plan of the article is the following. In Section \ref{Secc_ProofTh1} we give the proof of the Theorem \ref{Theoreme1}, in Section \ref{Secc_ProofTh2}, we study Theorem \ref{Theoreme2} and Section \ref{Secc_ProofTh3} is devoted to the proof of the Theorem \ref{Theoreme3}. Finally, in the Appendix \ref{Secc_ProofProposition1} we give the proof of the Proposition \ref{Proposition1} which gives an estimate for the operator $T^\alpha$ defined in the expression (\ref{Def_OperatorAlpha}). 
%%%%%%%%%%%%%%%%%%%%%%%%%%%%%%%%%%%%%%%%%%%%%%%%%%%
\mysection{Proof of the Theorem \ref{Theoreme1}}\label{Secc_ProofTh1}
We first consider the integral formulation of the equation (\ref{EquationIntro}) and we have 
\begin{equation}\label{Formulation_Integrale}
\theta(t,x)=\mathfrak{g}_t\ast \theta_0(x)+\int_{0}^{t}\mathfrak{g}_{t-s}\ast (\mathbb{T}_{[\theta]}\cdot \nabla \theta)(s,x)ds+\int_{0}^{t}\mathfrak{g}_{t-s}\ast f(s,x)ds.
\end{equation}
Applying the norm $\|\cdot\|_{E}$ defined in the formula (\ref{Espace_Resolution}) above to the previous expression we have 
\begin{equation}\label{Estimations_PointFixe12}
\|\theta\|_{E}\leq \underbrace{\|\mathfrak{g}_t\ast \theta_0\|_{E}}_{(1)}+\underbrace{\left\|\int_{0}^{t}\mathfrak{g}_{t-s}\ast (\mathbb{T}_{[\theta]}\cdot \nabla \theta)(s,\cdot)ds\right\|_{E}}_{(2)}+\underbrace{\left\|\int_{0}^{t}\mathfrak{g}_{t-s}\ast f(s,\cdot)ds\right\|_{E}}_{(3)},
\end{equation}
and we will study each term separately. 
\begin{itemize}
\item[(1)] The study of the initial data $\theta_0$ in (\ref{Estimations_PointFixe12}) is straightforward as we have by the hypothesis (\ref{Espace_ConditionInitiale}) the estimate
\begin{equation}\label{Estimation_DonneeInitiale}
\|\mathfrak{g}_t\ast \theta_0\|_{E}=\underset{t>0}{\sup}\; t^{\frac{2q-n}{2q}}\|\nabla (\mathfrak{g}_t\ast \theta_0)\|_{L^{q}}=\underset{t>0}{\sup}\; t^{\frac{2q-n}{2q}}\|\mathfrak{g}_t\ast \theta_0\|_{\dot{W}^{1,q}}=\|\theta_0\|_{\dot{B}^{\frac{n-q}{q},q}_\infty}<+\infty.
\end{equation}
%%%%%%%%%%%%%%%%%%%%%%%%%%%%%%%%%%%%%%%%%%%%%%%%%%%
\item[(2)] For the nonlinear drift given in $(2)$ of the expression (\ref{Estimations_PointFixe12}) above, we have: 
\begin{eqnarray*}
\left\|\int_{0}^{t}\mathfrak{g}_{t-s}\ast (\mathbb{T}_{[\theta]}\cdot \nabla \theta)(s,\cdot)ds\right\|_{E}&=&\underset{t>0}{\sup}\; t^{\frac{2q-n}{2q}}\left\|\nabla \int_{0}^{t}\mathfrak{g}_{t-s}\ast (\mathbb{T}_{[\theta]}\cdot \nabla \theta)(s,\cdot)ds\right\|_{L^{q}}\\
&=& \underset{t>0}{\sup}\; t^{\frac{2q-n}{2q}}\left\| \int_{0}^{t}\nabla\mathfrak{g}_{t-s}\ast (\mathbb{T}_{[\theta]}\cdot \nabla \theta)(s,\cdot)ds\right\|_{L^{q}},
\end{eqnarray*}
and we obtain
\begin{eqnarray*}
\left\|\int_{0}^{t}\mathfrak{g}_{t-s}\ast (\mathbb{T}_{[\theta]}\cdot \nabla \theta)(s,\cdot)ds\right\|_{E}&\leq&\underset{t>0}{\sup}\; t^{\frac{2q-n}{2q}} \int_{0}^{t}\|\nabla\mathfrak{g}_{t-s}\|_{L^r}\| (\mathbb{T}_{[\theta]}\cdot \nabla \theta)(s,\cdot)\|_{L^{p}}ds,
\end{eqnarray*}
where we applied the Young inequalities for the convolution with $1+\frac{1}{q}=\frac{1}{r}+\frac{1}{p}$.
Recall that by the usual estimates for the heat kernel we have 
$$\|\nabla\mathfrak{g}_{t-s}\|_{L^r}\leq C(t-s)^{-\frac{1+n(1-\frac{1}{r})}{2}}=C(t-s)^{-\frac{1+n(\frac{1}{p}-\frac{1}{q})}{2}}.$$
We thus obtain the estimate
$$\left\|\int_{0}^{t}\mathfrak{g}_{t-s}\ast (\mathbb{T}_{[\theta]}\cdot \nabla \theta)(s,\cdot)ds\right\|_{E}\leq C \,\underset{t>0}{\sup}\; t^{\frac{2q-n}{2q}} \int_{0}^{t} (t-s)^{-\frac{1+n(\frac{1}{p}-\frac{1}{q})}{2}}\| (\mathbb{T}_{[\theta]}\cdot \nabla \theta)(s,\cdot)\|_{L^{p}}ds,$$
and by the H\"older inequality with $\frac{1}{p}=\frac{1}{q}+\frac{1}{\mathfrak{s}}$, we can write
$$\left\|\int_{0}^{t}\mathfrak{g}_{t-s}\ast (\mathbb{T}_{[\theta]}\cdot \nabla \theta)(s,\cdot)ds\right\|_{E}\leq C \,\underset{t>0}{\sup}\; t^{\frac{2q-n}{2q}} \int_{0}^{t} (t-s)^{-\frac{1+n(\frac{1}{p}-\frac{1}{q})}{2}}\|\mathbb{T}_{[\theta]}(s,\cdot)\|_{L^{\mathfrak{s}}} \|\nabla \theta(s,\cdot)\|_{L^{q}}ds.$$
At this point we use the following key result 
%%%%%%%%%%%%%%%%%%%%%%%%%%%%%%%%%%%%%%%%%%%%%%%%%%%
\begin{Lemme}\label{LemmaRough} If $T_k$ is a rough singular integral operator of the form given in (\ref{Def_Operator}) associated to kernels $\Omega_k \in L^\rho(\mathbb{S}^{n-1})$ with $1<\rho<n$. If $\psi:\mathbb{R}^n\longrightarrow \mathbb{R}$ is a function such that $\nabla \psi\in L^q(\mathbb{R}^n)$ with $1<\frac{\rho n}{\rho n+\rho-n}\leq q<n$, then for all $1\leq k\leq n$ we have the estimate
$$\|T_k(\psi)\|_{L^\mathfrak{s}}\leq C\|\nabla\psi\|_{L^q},$$
where $\mathfrak{s}=\frac{nq}{n-q}>1$.
\end{Lemme}
See a proof of this inequality in \cite{ChMarcociMarcoci}. 
\begin{Remarque}
This functional estimate is the key to close the fixed point argument. To the best of our knowledge, we do not know if a better estimate is available for this type of rough operators $T_k$ in the context of Lebesgue spaces. 
\end{Remarque}
%%%%%%%%%%%%%%%%%%%%%%%%%%%%%%%%%%%%%%%%%%%%%%%%%%%
Thus, applying this lema we obtain the inequality
$$\left\|\int_{0}^{t}\mathfrak{g}_{t-s}\ast (\mathbb{T}_{[\theta]}\cdot \nabla \theta)(s,\cdot)ds\right\|_{E}\leq C \,\underset{t>0}{\sup}\; t^{\frac{2q-n}{2q}} \int_{0}^{t} (t-s)^{-\frac{1+n(\frac{1}{p}-\frac{1}{q})}{2}}\|\nabla\theta(s,\cdot)\|_{L^q}\|\nabla \theta(s,\cdot)\|_{L^{q}}ds,$$
which we rewrite as follows
\begin{eqnarray*}
\left\|\int_{0}^{t}\mathfrak{g}_{t-s}\ast (\mathbb{T}_{[\theta]}\cdot \nabla \theta)(s,\cdot)ds\right\|_{E}&\leq &C \,\underset{t>0}{\sup}\; t^{\frac{2q-n}{2q}} \int_{0}^{t} (t-s)^{-\frac{1+n(\frac{1}{p}-\frac{1}{q})}{2}}s^{-\frac{2q-n}{q}}\times \\
&&\left(s^{\frac{2q-n}{2q}}\|\nabla\theta(s,\cdot)\|_{L^q}\right)\left(s^{\frac{2q-n}{2q}}\|\nabla\theta(s,\cdot)\|_{L^q}\right)ds,
\end{eqnarray*}
and we have
\begin{equation*}
\begin{split}
\left\|\int_{0}^{t}\mathfrak{g}_{t-s}\ast (\mathbb{T}_{[\theta]}\cdot \nabla \theta)(s,\cdot)ds\right\|_{E}\leq C \,\underset{s>0}{\sup}\;\left(s^{\frac{2q-n}{2q}}\|\nabla\theta(s,\cdot)\|_{L^q}\right)\times \,\underset{s>0}{\sup}\;\left(s^{\frac{2q-n}{2q}}\|\nabla\theta(s,\cdot)\|_{L^q}\right) \\
\times \,\underset{t>0}{\sup}\; t^{\frac{2q-n}{2q}} \int_{0}^{t} (t-s)^{-\frac{1+n(\frac{1}{p}-\frac{1}{q})}{2}}s^{-\frac{2q-n}{q}} ds,
\end{split}
\end{equation*}
and using the definition of the norm of the space $E$ given in the expression (\ref{Espace_Resolution}) above, we obtain
$$\left\|\int_{0}^{t}\mathfrak{g}_{t-s}\ast (\mathbb{T}_{[\theta]}\cdot \nabla \theta)(s,\cdot)ds\right\|_{E}\leq C\|\theta\|_{E}\|\theta\|_{E}
\times \,\underset{t>0}{\sup}\; t^{\frac{2q-n}{2q}} \int_{0}^{t} (t-s)^{-\frac{1+n(\frac{1}{p}-\frac{1}{q})}{2}}s^{-\frac{2q-n}{q}} ds.$$
Recalling that we have $\frac{1}{p}=\frac{1}{q}+\frac{1}{\mathfrak{s}}$ and $\frac{1}{\mathfrak{s}}=\frac{1}{q}-\frac{1}{n}$, we deduce that $\frac{1}{p}-\frac{1}{q}=\frac{1}{q}-\frac{1}{n}$, so we can write
\begin{eqnarray*}
\left\|\int_{0}^{t}\mathfrak{g}_{t-s}\ast (\mathbb{T}_{[\theta]}\cdot \nabla \theta)(s,\cdot)ds\right\|_{E}&\leq &C\|\theta\|_{E}\|\theta\|_{E}
\times \,\underset{t>0}{\sup}\; t^{\frac{2q-n}{2q}} \int_{0}^{t} (t-s)^{-\frac{1+n(\frac{1}{q}-\frac{1}{n})}{2}}s^{-\frac{2q-n}{q}} ds\\
&\leq &C\|\theta\|_{E}\|\theta\|_{E}
\times \,\underset{t>0}{\sup}\; t^{\frac{2q-n}{2q}} \int_{0}^{t} (t-s)^{-\frac{n}{2q}}s^{-\frac{2q-n}{q}} ds.
\end{eqnarray*}
We study now the integral above and we have 
$$ \int_{0}^{t} (t-s)^{-\frac{n}{2q}}s^{-\frac{2q-n}{q}} ds= \int_{0}^{\frac{t}{2}} (t-s)^{-\frac{n}{2q}}s^{-\frac{2q-n}{q}} ds+ \int_{\frac{t}{2}}^{t} (t-s)^{-\frac{n}{2q}}s^{-\frac{2q-n}{q}}ds,$$
since if $0<s<\frac{t}{2}$, we have $\frac{t}{2}<(t-s)<t$ and $(t-s)^{-\frac{n}{2q}}<Ct^{-\frac{n}{2q}}$, and if $\frac{t}{2}<s<t$, we have $s^{-\frac{2q-n}{q}}<Ct^{-\frac{2q-n}{q}}$, we then obtain
$$ \int_{0}^{t} (t-s)^{-\frac{n}{2q}}s^{-\frac{2q-n}{q}} ds\leq C t^{-\frac{n}{2q}} \int_{0}^{\frac{t}{2}} s^{-\frac{2q-n}{q}} ds+C t^{-\frac{2q-n}{q}}\int_{\frac{t}{2}}^{t} (t-s)^{-\frac{n}{2q}}ds,$$
recalling that $1<q<n<2q$, the previous integrals are finite and we have
$$ \int_{0}^{t} (t-s)^{-\frac{n}{2q}}s^{-\frac{2q-n}{q}} ds\leq C t^{1-\frac{n}{2q}-\frac{2q-n}{q}}=Ct^{-\frac{2q-n}{2q}},$$
so we can write 
\begin{eqnarray}
\left\|\int_{0}^{t}\mathfrak{g}_{t-s}\ast (\mathbb{T}_{[\theta]}\cdot \nabla \theta)(s,\cdot)ds\right\|_{E}&\leq &C\|\theta\|_{E}\|\theta\|_{E}
\times \,\underset{t>0}{\sup}\; t^{\frac{2q-n}{2q}}\left( \int_{0}^{t} (t-s)^{-\frac{n}{2q}}s^{-\frac{2q-n}{2q}} ds\right)\notag\\
&\leq &  C\|\theta\|_{E}\|\theta\|_{E}\times \,\underset{t>0}{\sup}\; t^{\frac{2q-n}{2q}}\times t^{-\frac{2q-n}{2q}}\notag \\
&\leq & C_B\|\theta\|_{E}\|\theta\|_{E}, \label{Estimation_Bilineaire}
\end{eqnarray}
and this estimate gives the continuity of the nonlinear term in the space $E$.\\
%%%%%%%%%%%%%%%%%%%%%%%%%%%%%%%%%%%%%%%%%%%%%%%%%%%
\item[(3)] For the external force $f$ we write 
$$\left\|\int_{0}^{t}\mathfrak{g}_{t-s}\ast f(s,\cdot)ds\right\|_{E}=\underset{t>0}{\sup}\; t^{\frac{2q-n}{2q}}\left\|\nabla \int_{0}^{t}\mathfrak{g}_{t-s}\ast f(s,\cdot)ds \right\|_{L^{q}}\leq \underset{t>0}{\sup}\; t^{\frac{2q-n}{2q}} \int_{0}^{t}\|\nabla\mathfrak{g}_{t-s}\ast f(s,\cdot)\|_{L^q}ds,$$
and by the Young inequalities for the convolution with $1+\frac{1}{q}=\frac{1}{r}+\frac{1}{\varrho}$, we obtain
$$\left\|\int_{0}^{t}\mathfrak{g}_{t-s}\ast f(s,\cdot)ds\right\|_{E}\leq  \underset{t>0}{\sup}\; t^{\frac{2q-n}{2q}} \int_{0}^{t}\|\nabla\mathfrak{g}_{t-s}\|_{L^r}\|f(s,\cdot)\|_{L^\varrho}ds,$$
and since we have the control $\|\nabla\mathfrak{g}_{t-s}\|_{L^r}\leq C(t-s)^{-\frac{1+n(1-\frac{1}{r})}{2}}=C(t-s)^{-\frac{1+n(\frac{1}{\varrho}-\frac{1}{q})}{2}}$ for the heat kernel, we have 
$$\left\|\int_{0}^{t}\mathfrak{g}_{t-s}\ast f(s,\cdot)ds\right\|_{E}\leq  C\underset{t>0}{\sup}\; t^{\frac{2q-n}{2q}} \int_{0}^{t} (t-s)^{-\frac{1+n(\frac{1}{\varrho}-\frac{1}{q})}{2}}\|f(s,\cdot)\|_{L^\varrho}ds.$$
In order to obtain the wished norm $\|\cdot\|_{E_f}$ given in the expression (\ref{Espace_ForceExterieure}) for the external force we write 
\begin{eqnarray*}
\left\|\int_{0}^{t}\mathfrak{g}_{t-s}\ast f(s,\cdot)ds\right\|_{E}&\leq &C\underset{t>0}{\sup}\; t^{\frac{2q-n}{2q}} \int_{0}^{t} (t-s)^{-\frac{1+n(\frac{1}{\varrho}-\frac{1}{q})}{2}}s^{-(\frac{3}{2}-\frac{n}{2\varrho})}\left(s^{\frac{3}{2}-\frac{n}{2\varrho}}\|f(s,\cdot)\|_{L^\rho}\right)ds\\
&\leq & C \underset{s>0}{\sup}\left(s^{\frac{3}{2}-\frac{n}{2\varrho}}\|f(s,\cdot)\|_{L^\varrho}\right)\;\;\underset{t>0}{\sup}\; t^{\frac{2q-n}{2q}} \int_{0}^{t} (t-s)^{-\frac{1+n(\frac{1}{\varrho}-\frac{1}{q})}{2}}s^{-(\frac{3}{2}-\frac{n}{2\varrho})}ds\\
&\leq & C \|f\|_{E_f}\;\;\underset{t>0}{\sup}\; t^{\frac{2q-n}{2q}} \int_{0}^{t} (t-s)^{-\frac{1+n(\frac{1}{\varrho}-\frac{1}{q})}{2}}s^{-(\frac{3}{2}-\frac{n}{2\varrho})}ds.
\end{eqnarray*}
We now have 
$$\int_{0}^{t} (t-s)^{-\frac{1+n(\frac{1}{\varrho}-\frac{1}{q})}{2}}s^{-(\frac{3}{2}-\frac{n}{2\varrho})}ds=\int_{0}^{t/2} (t-s)^{-\frac{1+n(\frac{1}{\varrho}-\frac{1}{q})}{2}}s^{-(\frac{3}{2}-\frac{n}{2\rho})}ds+\int_{t/2}^{t} (t-s)^{-\frac{1+n(\frac{1}{\varrho}-\frac{1}{q})}{2}}s^{-(\frac{3}{2}-\frac{n}{2\rho})}ds,$$
but since, if $0<s<\frac{t}{2}$, we have $\frac{t}{2}<(t-s)<t$ and $(t-s)^{-\frac{1+n(\frac{1}{\rho}-\frac{1}{q})}{2}}<Ct^{-\frac{1+n(\frac{1}{\varrho}-\frac{1}{q})}{2}}$, and if $\frac{t}{2}<s<t$, we have $s^{-(\frac{3}{2}-\frac{n}{2\varrho})}<Ct^{-(\frac{3}{2}-\frac{n}{2\varrho})}$, we then obtain
$$\int_{0}^{t} (t-s)^{-\frac{1+n(\frac{1}{\varrho}-\frac{1}{q})}{2}}s^{-(\frac{3}{2}-\frac{n}{2\varrho})}ds\leq  t^{-\frac{1+n(\frac{1}{\varrho}-\frac{1}{q})}{2}}\int_{0}^{t/2}s^{-(\frac{3}{2}-\frac{n}{2\varrho})}ds+t^{-(\frac{3}{2}-\frac{n}{2\varrho})}\int_{t/2}^{t} (t-s)^{-\frac{1+n(\frac{1}{\varrho}-\frac{1}{q})}{2}}ds.$$
At this point we note that the hypotheses $\varrho<n<3\varrho$ and $\frac{1}{q}-\frac{1}{n}<\frac{1}{\varrho}<\frac{1}{q}$ guarantee that the two previous integrals are bounded and we can write
$$\int_{0}^{t} (t-s)^{-\frac{1+n(\frac{1}{\varrho}-\frac{1}{q})}{2}}s^{-(\frac{3}{2}-\frac{n}{2\varrho})}ds\leq  Ct^{-\frac{1+n(\frac{1}{\varrho}-\frac{1}{q})}{2}}t^{1-(\frac{3}{2}-\frac{n}{2\varrho})}=Ct^{-\frac{2q-n}{2q}}.$$
Now, we finally obtain 
\begin{eqnarray} 
\left\|\int_{0}^{t}\mathfrak{g}_{t-s}\ast f(s,\cdot)ds\right\|_{E}&\leq & C \|f\|_{E_f}\;\;\underset{t>0}{\sup}\; t^{\frac{2q-n}{2q}} \int_{0}^{t} (t-s)^{-\frac{1+n(\frac{1}{\varrho}-\frac{1}{q})}{2}}s^{-(\frac{3}{2}-\frac{n}{2\varrho})}ds\notag \\
&\leq & C \|f\|_{E_f}\;\;\underset{t>0}{\sup}\; t^{\frac{2q-n}{2q}} \times t^{-\frac{2q-n}{2q}}\leq C \|f\|_{E_f}<+\infty.\label{Estimation_ForceExterieure}
\end{eqnarray}
\end{itemize}
Now with the estimates (\ref{Estimation_DonneeInitiale}), (\ref{Estimation_ForceExterieure}) and (\ref{Estimation_Bilineaire}), as long as we have that the quantity  $\|\theta_0\|_{\dot{B}^{\frac{n-q}{q},q}_\infty}+\|f\|_{E_f}$ is small enough, we can apply the Banach-Picard fixed point argument to obtain a global mild solution of the integral problem (\ref{Formulation_Integrale}) and the Theorem \ref{Theoreme1} is now proven. \hfill $\blacksquare$
%%%%%%%%%%%%%%%%%%%%%%%%%%%%%%%%%%%%%%%%%%%%%%%%%%%
\mysection{Proof of the Theorem \ref{Theoreme2}}\label{Secc_ProofTh2}
We consider the integral formulation of the equation (\ref{EquationIntroDiv})
\begin{equation}\label{Formulation_IntegraleDiv}
\theta(t,x)=\mathfrak{g}_t\ast \theta_0(x)+\int_{0}^{t}\mathfrak{g}_{t-s}\ast div(\mathbb{T}_{[(-\Delta)^{-\frac{1}{2}}\theta]}\theta)(s,x)ds+\int_{0}^{t}\mathfrak{g}_{t-s}\ast f(s,x)ds,
\end{equation}
and applying the norm $\|\cdot\|_{\mathbb{E}}$ we obtain
\begin{equation}\label{Estimations_PointFixeDiv}
\|\theta\|_{\mathbb{E}}\leq \underbrace{\|\mathfrak{g}_t\ast \theta_0\|_{\mathbb{E}}}_{(1)}+\underbrace{\left\|\int_{0}^{t}\mathfrak{g}_{t-s}\ast div(\mathbb{T}_{[(-\Delta)^{-\frac{1}{2}}\theta]}\theta)(s,\cdot)ds\right\|_{\mathbb{E}}}_{(2)}+\underbrace{\left\|\int_{0}^{t}\mathfrak{g}_{t-s}\ast f(s,\cdot)ds\right\|_{\mathbb{E}}}_{(3)}.
\end{equation}
As before, we will study each term above separately. 
\begin{itemize}
\item For the initial data $\theta_0$ in (\ref{Estimations_PointFixeDiv}) we simply write
\begin{equation}\label{Estimation_DonneeInitialeDiv}
\|\mathfrak{g}_t\ast \theta_0\|_{\mathbb{E}}=\underset{t>0}{\sup}\; t^{\frac{2q-n}{2q}}\|\mathfrak{g}_t\ast \theta_0\|_{L^{q}}=\|\theta_0\|_{\dot{B}^{-\frac{2q-n}{q},q}_\infty},
\end{equation}
which a bounded quantity by the hypothesis (\ref{Espace_ConditionInitialeDiv}).
%%%%%%%%%%%%%%%%%%%%%%%%%%%%%%%%%%%%%%%%%%%%%%%%%%%
\item For the nonlinear drift given in $(2)$ of the expression (\ref{Estimations_PointFixeDiv}) above, we have: 
$$\left\|\int_{0}^{t}\mathfrak{g}_{t-s}\ast div(\mathbb{T}_{[(-\Delta)^{-\frac{1}{2}}\theta]}\theta)(s,\cdot)ds\right\|_{\mathbb{E}}=\underset{t>0}{\sup}\; t^{\frac{2q-n}{2q}}\left\| \int_{0}^{t}\mathfrak{g}_{t-s}\ast div(\mathbb{T}_{[(-\Delta)^{-\frac{1}{2}}\theta]}\theta)(s,\cdot)ds\right\|_{L^{q}},$$
and applying the Young inequalities with $1+\frac{1}{q}=\frac{1}{r}+\frac{1}{p}$, we get
$$\left\|\int_{0}^{t}\mathfrak{g}_{t-s}\ast div(\mathbb{T}_{[(-\Delta)^{-\frac{1}{2}}\theta]}\theta)(s,\cdot)ds \right\|_{\mathbb{E}} \leq  \underset{t>0}{\sup}\; t^{\frac{2q-n}{2q}} \int_{0}^{t}\|\nabla\mathfrak{g}_{t-s}\|_{L^r}\|(\mathbb{T}_{[(-\Delta)^{-\frac{1}{2}}\theta]}\theta)(s,\cdot)\|_{L^{p}}ds.$$
Since we have the control $\|\nabla\mathfrak{g}_{t-s}\|_{L^r}\leq C(t-s)^{-\frac{1+n(1-\frac{1}{r})}{2}}=C(t-s)^{-\frac{1+n(\frac{1}{p}-\frac{1}{q})}{2}}$, we can write
$$\left\|\int_{0}^{t}\mathfrak{g}_{t-s}\ast div(\mathbb{T}_{[(-\Delta)^{-\frac{1}{2}}\theta]}\theta)(s,\cdot)ds\right\|_{\mathbb{E}}\leq C \,\underset{t>0}{\sup}\; t^{\frac{2q-n}{2q}} \int_{0}^{t} (t-s)^{-\frac{1+n(\frac{1}{p}-\frac{1}{q})}{2}}\| (\mathbb{T}_{[(-\Delta)^{-\frac{1}{2}}\theta]}\theta)(s,\cdot)\|_{L^{p}}ds,$$
and by the H\"older inequality with $\frac{1}{p}=\frac{1}{q}+\frac{1}{\mathfrak{s}}$, we can write
$$\left\|\int_{0}^{t}\mathfrak{g}_{t-s}\ast div(\mathbb{T}_{[(-\Delta)^{-\frac{1}{2}}\theta]}\theta)(s,\cdot)ds\right\|_{\mathbb{E}}\leq C \,\underset{t>0}{\sup}\; t^{\frac{2q-n}{2q}} \int_{0}^{t} (t-s)^{-\frac{1+n(\frac{1}{p}-\frac{1}{q})}{2}}\|\mathbb{T}_{[(-\Delta)^{-\frac{1}{2}}\theta]}\|_{L^{\mathfrak{s}}} \|\theta(s,\cdot)\|_{L^{q}}ds.$$
At this point we use the Lemma \ref{LemmaRough}, which gives to us the control 
$$\|\mathbb{T}_{[(-\Delta)^{-\frac{1}{2}}\theta]}\|_{L^{\mathfrak{s}}} \leq C\|\nabla((-\Delta)^{-\frac{1}{2}}\theta)(s,\cdot)\|_{L^q},$$
where $\frac{1}{\mathfrak{s}}=\frac{1}{q}-\frac{1}{n}$. Since $1<q<n$ and since the Riesz transforms are bounded in the Lebesgue spaces $L^q(\mathbb{R}^n)$, we obtain the inequality
$$\left\|\int_{0}^{t}\mathfrak{g}_{t-s}\ast div(\mathbb{T}_{[(-\Delta)^{-\frac{1}{2}}\theta]}\theta)(s,\cdot)ds\right\|_{\mathbb{E}}\leq C \,\underset{t>0}{\sup}\; t^{\frac{2q-n}{2q}} \int_{0}^{t} (t-s)^{-\frac{1+n(\frac{1}{p}-\frac{1}{q})}{2}}\|\theta(s,\cdot)\|_{L^q}\| \theta(s,\cdot)\|_{L^{q}}ds,$$
which we rewrite as follows
\begin{eqnarray*}
\left\|\int_{0}^{t}\mathfrak{g}_{t-s}\ast div(\mathbb{T}_{[(-\Delta)^{-\frac{1}{2}}\theta]}\theta)(s,\cdot)ds\right\|_{\mathbb{E}}&\leq &C \,\underset{t>0}{\sup}\; t^{\frac{2q-n}{2q}} \int_{0}^{t} (t-s)^{-\frac{1+n(\frac{1}{p}-\frac{1}{q})}{2}}s^{-\frac{2q-n}{q}}\times \\
&&\left(s^{\frac{2q-n}{2q}}\|\theta(s,\cdot)\|_{L^q}\right)\left(s^{\frac{2q-n}{2q}}\|\theta(s,\cdot)\|_{L^q}\right)ds,
\end{eqnarray*}
and we have
\begin{equation*}
\begin{split}
\left\|\int_{0}^{t}\mathfrak{g}_{t-s}\ast div(\mathbb{T}_{[(-\Delta)^{-\frac{1}{2}}\theta]}\theta)(s,\cdot)ds\right\|_{\mathbb{E}}\leq C \,\underset{s>0}{\sup}\;\left(s^{\frac{2q-n}{2q}}\|\theta(s,\cdot)\|_{L^q}\right)\times \,\underset{s>0}{\sup}\;\left(s^{\frac{2q-n}{2q}}\|\theta(s,\cdot)\|_{L^q}\right) \\
\times \,\underset{t>0}{\sup}\; t^{\frac{2q-n}{2q}} \int_{0}^{t} (t-s)^{-\frac{1+n(\frac{1}{p}-\frac{1}{q})}{2}}s^{-\frac{2q-n}{q}} ds.
\end{split}
\end{equation*}
Now using the definition of the norm of the space $\mathbb{E}$ given in the expression (\ref{Espace_ResolutionDiv}) above, we obtain
$$\left\|\int_{0}^{t}\mathfrak{g}_{t-s}\ast div(\mathbb{T}_{[(-\Delta)^{-\frac{1}{2}}\theta]}\theta)(s,\cdot)ds\right\|_{\mathbb{E}}\leq C\|\theta\|_{\mathbb{E}}\|\theta\|_{\mathbb{E}}
\times \,\underset{t>0}{\sup}\; t^{\frac{2q-n}{2q}} \int_{0}^{t} (t-s)^{-\frac{1+n(\frac{1}{p}-\frac{1}{q})}{2}}s^{-\frac{2q-n}{q}} ds.$$
By the same computation performed in the previous theorem, we have that 
$$\underset{t>0}{\sup}\; t^{\frac{2q-n}{2q}} \int_{0}^{t} (t-s)^{-\frac{1+n(\frac{1}{p}-\frac{1}{q})}{2}}s^{-\frac{2q-n}{q}} ds<+\infty,$$ 
so we can write 
\begin{equation}\label{Estimation_BilineaireDiv}
\left\|\int_{0}^{t}\mathfrak{g}_{t-s}\ast div(\mathbb{T}_{[(-\Delta)^{-\frac{1}{2}}\theta]}\theta)(s,\cdot)ds\right\|_{\mathbb{E}}\leq  C_B\|\theta\|_{\mathbb{E}}\|\theta\|_{\mathbb{E}}, 
\end{equation}
and this estimate gives the continuity of the nonlinear term in the space $\mathbb{E}$.
%%%%%%%%%%%%%%%%%%%%%%%%%%%%%%%%%%%%%%%%%%%%%%%%%%%
\item The last term (3) of the expression (\ref{Estimations_PointFixeDiv}) is treated as follows:
\begin{eqnarray*}
\left\|\int_{0}^{t}\mathfrak{g}_{t-s}\ast f(s,\cdot)ds\right\|_{\mathbb{E}}&=&\underset{t>0}{\sup}\; t^{\frac{2q-n}{2q}}\left\| \int_{0}^{t}\mathfrak{g}_{t-s}\ast f(s,\cdot)ds\right\|_{L^{q}}\\
&=&\underset{t>0}{\sup}\; t^{\frac{2q-n}{2q}}\left\| \int_{0}^{t}\mathfrak{g}_{t-s}\ast (-\Delta)^{\frac{1}{2}}(-\Delta)^{-\frac{1}{2}} f(s,\cdot)ds\right\|_{L^{q}},
\end{eqnarray*}
and by the Young inequalities with $1+\frac{1}{q}=\frac{1}{r}+\frac{1}{\varrho}$, we have
\begin{eqnarray*}
\left\|\int_{0}^{t}\mathfrak{g}_{t-s}\ast f(s,\cdot)ds\right\|_{\mathbb{E}}&\leq &\underset{t>0}{\sup}\; t^{\frac{2q-n}{2q}} \int_{0}^{t}\|(-\Delta)^{\frac{1}{2}}\mathfrak{g}_{t-s}\|_{L^r} \|(-\Delta)^{-\frac{1}{2}} f(s,\cdot)\|_{L^\varrho}ds\\
&\leq & \underset{t>0}{\sup}\; t^{\frac{2q-n}{2q}} \int_{0}^{t}(t-s)^{-\frac{1+n(\frac{1}{\varrho}-\frac{1}{q})}{2}}\|f(s,\cdot)\|_{\dot{W}^{-1,\varrho}}ds,
\end{eqnarray*}
where we used the usual estimates for the heat kernel. In order to reconstruct the norm $\|\cdot\|_{\mathbb{E}_f}$ given in the formula (\ref{Espace_ForceExterieureDiv}), we write
\begin{eqnarray*}
\left\|\int_{0}^{t}\mathfrak{g}_{t-s}\ast f(s,\cdot)ds\right\|_{\mathbb{E}}&\leq & \underset{t>0}{\sup}\; t^{\frac{2q-n}{2q}} \int_{0}^{t}(t-s)^{-\frac{1+n(\frac{1}{\varrho}-\frac{1}{q})}{2}}s^{-(\frac{3}{2}-\frac{n}{2\varrho})} \left(s^{(\frac{3}{2}-\frac{n}{2\varrho})} \|f(s,\cdot)\|_{\dot{W}^{-1,\varrho}}\right)ds\\
&\leq & \underset{s>0}{\sup}\;\left(s^{(\frac{3}{2}-\frac{n}{2\varrho})} \|f(s,\cdot)\|_{\dot{W}^{-1,\varrho}}\right)\; \underset{t>0}{\sup}\; t^{\frac{2q-n}{2q}} \int_{0}^{t}(t-s)^{-\frac{1+n(\frac{1}{\varrho}-\frac{1}{q})}{2}}s^{-(\frac{3}{2}-\frac{n}{2\varrho})} ds\\
&\leq & \|f\|_{\mathbb{E}_f}\; \underset{t>0}{\sup}\; t^{\frac{2q-n}{2q}} \int_{0}^{t}(t-s)^{-\frac{1+n(\frac{1}{\varrho}-\frac{1}{q})}{2}}s^{-(\frac{3}{2}-\frac{n}{2\varrho})} ds.
\end{eqnarray*}
Since by hypothesis we have $1<\varrho<n<3\varrho$ and $\frac{1}{q}-\frac{1}{n}<\frac{1}{\varrho}<\frac{1}{q}$ then, by the same computations as before we have $\displaystyle{\underset{t>0}{\sup}\; t^{\frac{2q-n}{2q}} \int_{0}^{t}(t-s)^{-\frac{1+n(\frac{1}{\varrho}-\frac{1}{q})}{2}}s^{-\frac{2q-n}{q}} ds}<+\infty$, and we can write
\begin{equation}\label{Estimation_ForceExterieureDiv}
\left\|\int_{0}^{t}\mathfrak{g}_{t-s}\ast f(s,\cdot)ds\right\|_{\mathbb{E}}\leq C\|f\|_{\mathbb{E}_f}.
\end{equation}
\end{itemize}
%%%%%%%%%%%%%%%%%%%%%%%%%%%%%%%%%%%%%%%%%%%%%%%%%%%
With the estimates (\ref{Estimation_DonneeInitialeDiv}), (\ref{Estimation_BilineaireDiv}) and (\ref{Estimation_ForceExterieureDiv}) if the quantity $\|\theta_0\|_{\dot{B}^{-\frac{2q-n}{q},q}_\infty}+\|f\|_{\mathbb{E}_f}$ is small enough, by applying a fixed point argument we obtain a global mild solution of the integral problem (\ref{Formulation_IntegraleDiv}) and the Theorem \ref{Theoreme2} is now proven. \hfill $\blacksquare$
%%%%%%%%%%%%%%%%%%%%%%%%%%%%%%%%%%%%%%%%%%%%%%%%%%%
\mysection{Proof of the Theorem \ref{Theoreme3}}\label{Secc_ProofTh3}
We consider the following integral problem
\begin{equation}\label{Formulation_IntegraleAlpha}
\theta(t,x)=\mathfrak{g}_t\ast \theta_0(x)+\int_{0}^{t}\mathfrak{g}_{t-s}\ast (\mathbb{T}_{[\theta]}^{\alpha}\cdot \nabla \theta)(s,x)ds+\int_{0}^{t}\mathfrak{g}_{t-s}\ast f(s,x)ds,
\end{equation}
and applying the norm $\|\cdot\|_{\mathcal{E}}$ we obtain
\begin{equation}\label{Estimations_PointFixeAlpha}
\|\theta\|_{\mathcal{E}}\leq \underbrace{\|\mathfrak{g}_t\ast \theta_0\|_{\mathcal{E}}}_{(1)}+\underbrace{\left\|\int_{0}^{t}\mathfrak{g}_{t-s}\ast (\mathbb{T}_{[\theta]}^{\alpha}\cdot \nabla \theta)(s,\cdot)ds\right\|_{\mathcal{E}}}_{(2)}+\underbrace{\left\|\int_{0}^{t}\mathfrak{g}_{t-s}\ast f(s,\cdot)ds\right\|_{\mathcal{E}}}_{(3)}.
\end{equation}
As before, we will study each term above separately. 
\begin{itemize}
\item For the initial data $\theta_0$ in (\ref{Estimations_PointFixeAlpha}) we simply write
\begin{equation}\label{Estimation_DonneeInitialeAlpha}
\|\mathfrak{g}_t\ast \theta_0\|_{\mathcal{E}}=\underset{t>0}{\sup}\; t^{\frac{(2+\alpha)q-n}{2q}}\|\mathfrak{g}_{t}\ast \theta_0\|_{\dot{W}^{1,q}}=\|\theta_0\|_{\dot{B}^{\frac{n-(1+\alpha)q}{q},q}_\infty},
\end{equation}
which a bounded quantity by the hypothesis (\ref{Espace_ConditionInitialeAlpha}).
%%%%%%%%%%%%%%%%%%%%%%%%%%%%%%%%%%%%%%%%%%%%%%%%%%%
\item For the nonlinear drift given in $(2)$ we have: 
$$\left\|\int_{0}^{t}\mathfrak{g}_{t-s}\ast (\mathbb{T}_{[\theta]}^{\alpha}\cdot \nabla \theta)(s,\cdot)ds\right\|_{\mathcal{E}}=\underset{t>0}{\sup}\; t^{\frac{(2+\alpha)q-n}{2q}}\left\| \nabla\int_{0}^{t}\mathfrak{g}_{t-s}\ast(\mathbb{T}_{[\theta]}^{\alpha}\cdot \nabla \theta)(s,\cdot)ds\right\|_{L^{q}},$$
and applying the Young inequalities with $1+\frac{1}{q}=\frac{1}{r}+\frac{1}{p}$, we get
$$\left\|\int_{0}^{t}\mathfrak{g}_{t-s}\ast (\mathbb{T}_{[\theta]}^{\alpha}\cdot \nabla \theta)(s,\cdot)ds \right\|_{\mathcal{E}} \leq  \underset{t>0}{\sup}\;  t^{\frac{(2+\alpha)q-n}{2q}} \int_{0}^{t}\|\nabla\mathfrak{g}_{t-s}\|_{L^r}\|(\mathbb{T}_{[\theta]}^{\alpha}\nabla \theta)(s,\cdot)\|_{L^{p}}ds.$$
Since we have the control $\|\nabla\mathfrak{g}_{t-s}\|_{L^r}\leq C(t-s)^{-\frac{1+n(1-\frac{1}{r})}{2}}=C(t-s)^{-\frac{1+n(\frac{1}{p}-\frac{1}{q})}{2}}$, we can write
$$\left\|\int_{0}^{t}\mathfrak{g}_{t-s}\ast (\mathbb{T}_{[\theta]}^{\alpha}\cdot \nabla \theta)(s,\cdot)ds\right\|_{\mathcal{E}}\leq C \,\underset{t>0}{\sup}\;  t^{\frac{(2+\alpha)q-n}{2q}}\int_{0}^{t} (t-s)^{-\frac{1+n(\frac{1}{p}-\frac{1}{q})}{2}}\| (\mathbb{T}_{[\theta]}^{\alpha}\cdot \nabla \theta)(s,\cdot)\|_{L^{p}}ds,$$
and by the H\"older inequality with $\frac{1}{p}=\frac{1}{q}+\frac{1}{\mathfrak{s}}$, we can write
$$\left\|\int_{0}^{t}\mathfrak{g}_{t-s}\ast (\mathbb{T}_{[\theta]}^{\alpha}\cdot \nabla \theta)(s,\cdot)ds\right\|_{\mathcal{E}}\leq C \,\underset{t>0}{\sup}\;  t^{\frac{(2+\alpha)q-n}{2q}} \int_{0}^{t} (t-s)^{-\frac{1+n(\frac{1}{p}-\frac{1}{q})}{2}}\|\mathbb{T}_{[\theta]}^{\alpha}(s, \cdot)\|_{L^{\mathfrak{s}}} \|\nabla \theta(s,\cdot)\|_{L^{q}}ds.$$
We will use now the following result:
%%%%%%%%%%%%%%%%%%%%%%%%%%%%%%%%%%%%%%%%%%%%%%%%%%%
\begin{Proposition}\label{Proposition1} Consider a regular function $\phi:\mathbb{R}^n\longrightarrow \mathbb{R}$ and for some $0<\alpha<n-1$, consider the operator $T^\alpha_k$ as defined in the expression (\ref{Def_OperatorAlpha}) associated to a kernel $\Omega_k\in L^{\rho}(\mathbb{S}^{n-1})$ with $1<\rho<n$. Then we have the control 
$$\|T^\alpha_k(\phi)\|_{L^{\mathfrak{s}}}\leq C\|\nabla \phi\|_{L^q},$$
where $\mathfrak{s}=\frac{qn}{n-q(1+\alpha)}>1$ and $1<\frac{\rho n}{\rho n+\rho-n}<q<\frac{n}{1+\alpha}$.
\end{Proposition}
%%%%%%%%%%%%%%%%%%%%%%%%%%%%%%%%%%%%%%%%%%%%%%%%%%%
The proof of this result, interesting on its own, is postponed to the Appendix  \ref{Secc_ProofProposition1} below.\\

With this crucial estimate at hand, we can now write 
$$\left\|\int_{0}^{t}\mathfrak{g}_{t-s}\ast (\mathbb{T}_{[\theta]}^{\alpha}\cdot \nabla \theta)(s,\cdot)ds\right\|_{\mathcal{E}}\leq C\,\underset{t>0}{\sup}\; t^{\frac{(2+\alpha)q-n}{2q}}\int_{0}^{t} (t-s)^{-\frac{1+n(\frac{1}{p}-\frac{1}{q})}{2}}\|\nabla \theta(s, \cdot)\|_{L^{q}}\|\nabla\theta(s,\cdot)\|_{L^{q}}ds.$$
Noting that $\frac{1}{\mathfrak{s}}=\frac{1}{q}-\frac{1+\alpha}{n}$, due to the relationship $\frac{1}{p}=\frac{1}{q}+\frac{1}{\mathfrak{s}}$, we have that $\frac{1}{p}-\frac{1}{q}=\frac{1}{q}-\frac{1+\alpha}{n}$, and thus we obtain the identity $1+n(\frac{1}{p}-\frac{1}{q})=\frac{n-\alpha q}{q}$ which is a positive quantity since $1<\frac{\rho n}{\rho n+\rho-n}<q<\frac{n}{1+\alpha}$. Now, with this identity, we have 
$$\left\|\int_{0}^{t}\mathfrak{g}_{t-s}\ast (\mathbb{T}_{[\theta]}^{\alpha}\cdot \nabla \theta)(s,\cdot)ds\right\|_{\mathcal{E}}\leq C\,\underset{t>0}{\sup}\; t^{\frac{(2+\alpha)q-n}{2q}}\int_{0}^{t} (t-s)^{-(\frac{n-\alpha q}{2q})}\|\nabla \theta(s, \cdot)\|_{L^{q}}\|\nabla\theta(s,\cdot)\|_{L^{q}}ds,$$
which we rewrite as follows
\begin{eqnarray*}
\left\|\int_{0}^{t}\mathfrak{g}_{t-s}\ast (\mathbb{T}_{[\theta]}^{\alpha}\cdot \nabla \theta)(s,\cdot)ds\right\|_{\mathcal{E}}&\leq &C \,\underset{t>0}{\sup}\; t^{\frac{(2+\alpha)q-n}{2q}} \int_{0}^{t} (t-s)^{-(\frac{n-\alpha q}{2q})}s^{-(\frac{(2+\alpha)q-n}{q})}\times \\
&&\times \left(s^{\frac{(2+\alpha)q-n}{2q}}\|\nabla\theta(s,\cdot)\|_{L^q}\right)\left(s^{\frac{(2+\alpha)q-n}{2q}}\|\nabla\theta(s,\cdot)\|_{L^q}\right)ds.
\end{eqnarray*}
By the definition of the norm of the space $\mathcal{E}$ given in the expression (\ref{Espace_ResolutionAlpha}) above, we obtain
\begin{eqnarray*}
\left\|\int_{0}^{t}\mathfrak{g}_{t-s}\ast (\mathbb{T}_{[\theta]}^{\alpha}\cdot \nabla \theta)(s,\cdot)ds\right\|_{\mathcal{E}}&\leq &C\, \underset{s>0}{\sup}\;\left(s^{\frac{(2+\alpha)q-n}{2q}}\|\nabla\theta(s,\cdot)\|_{L^q}\right)\underset{s>0}{\sup}\;\left(s^{\frac{(2+\alpha)q-n}{2q}}\|\nabla\theta(s,\cdot)\|_{L^q}\right)\\
&&\,\underset{t>0}{\sup}\; t^{\frac{(2+\alpha)q-n}{2q}} \int_{0}^{t} (t-s)^{-(\frac{n-\alpha q}{2q})}s^{-(\frac{(2+\alpha)q-n}{q})}ds\\
&\leq &C\, \|\theta\|_{\mathcal{E}}\|\theta\|_{\mathcal{E}}\,\;\underset{t>0}{\sup}\; t^{\frac{(2+\alpha)q-n}{2q}} \int_{0}^{t} (t-s)^{-(\frac{n-\alpha q}{2q})}s^{-(\frac{(2+\alpha)q-n}{q})}ds.
\end{eqnarray*}
Recalling that $1<\max\{\frac{\rho n}{\rho n+\rho-n}, \frac{n}{2+\alpha}\}<q<\frac{n}{1+\alpha}$, by the same computation performed in the previous theorems, we have that $\displaystyle{\underset{t>0}{\sup}\; t^{\frac{(2+\alpha)q-n}{2q}} \int_{0}^{t} (t-s)^{-(\frac{n-\alpha q}{2q})}s^{-(\frac{(2+\alpha)q-n}{q})}ds<+\infty}$, so we can write 
\begin{equation}\label{Estimation_BilineaireAlpha}
\left\|\int_{0}^{t}\mathfrak{g}_{t-s}\ast (\mathbb{T}_{[\theta]}^{\alpha}\cdot \nabla \theta)(s,\cdot)ds\right\|_{\mathcal{E}}\leq  C_B\|\theta\|_{\mathcal{E}}\|\theta\|_{\mathcal{E}}, 
\end{equation}
and this estimate gives the continuity of the nonlinear term in the space $\mathcal{E}$.

%%%%%%%%%%%%%%%%%%%%%%%%%%%%%%%%%%%%%%%%%%%%%%%%%%%
\item For the external force, \emph{i.e.} the term (3) of the expression (\ref{Estimations_PointFixeAlpha}), we have
$$\left\|\int_{0}^{t}\mathfrak{g}_{t-s}\ast f(s,\cdot)ds\right\|_{\mathcal{E}}=\underset{t>0}{\sup}\; t^{\frac{(2+\alpha)q-n}{2q}}\left\| \nabla\int_{0}^{t}\mathfrak{g}_{t-s}\ast f(s,\cdot)ds\right\|_{L^{q}},$$
as before, by the Young inequalities with $1+\frac{1}{q}=\frac{1}{r}+\frac{1}{\varrho}$, we obtain by the usual estimates for the heat kernel:
\begin{eqnarray*}
\left\|\int_{0}^{t}\mathfrak{g}_{t-s}\ast f(s,\cdot)ds\right\|_{\mathcal{E}}&\leq &\underset{t>0}{\sup}\; t^{\frac{(2+\alpha)q-n}{2q}} \int_{0}^{t}\|\nabla\mathfrak{g}_{t-s}\|_{L^r} \|f(s,\cdot)\|_{L^\varrho}ds\\
&\leq & \underset{t>0}{\sup}\; t^{\frac{(2+\alpha)q-n}{2q}} \int_{0}^{t}(t-s)^{-\frac{1+n(\frac{1}{\varrho}-\frac{1}{q})}{2}}\|f(s,\cdot)\|_{L^{\varrho}}ds.
\end{eqnarray*}
In order to obtain the norm $\|\cdot\|_{\mathcal{E}_f}$ given in the formula (\ref{Espace_ForceExterieureAlpha}), we write
\begin{eqnarray*}
\left\|\int_{0}^{t}\mathfrak{g}_{t-s}\ast f(s,\cdot)ds\right\|_{\mathcal{E}}&\leq & \underset{t>0}{\sup}\; t^{\frac{(2+\alpha)q-n}{2q}} \int_{0}^{t}(t-s)^{-\frac{1+n(\frac{1}{\varrho}-\frac{1}{q})}{2}}s^{-(\frac{3+\alpha}{2}-\frac{n}{2\varrho})} \left(s^{(\frac{3+\alpha}{2}-\frac{n}{2\varrho})} \|f(s,\cdot)\|_{L^{\rho}}\right)ds\\
&\leq & \underset{s>0}{\sup}\;\left(s^{(\frac{3+\alpha}{2}-\frac{n}{2\varrho})} \|f(s,\cdot)\|_{L^{\varrho}}\right)\; \underset{t>0}{\sup}\; t^{\frac{(2+\alpha)q-n}{2q}} \int_{0}^{t}(t-s)^{-\frac{1+n(\frac{1}{\varrho}-\frac{1}{q})}{2}}s^{-(\frac{3+\alpha}{2}-\frac{n}{2\varrho})} ds\\
&\leq & \|f\|_{\mathcal{E}_f}\; \underset{t>0}{\sup}\; t^{\frac{(2+\alpha)q-n}{2q}} \int_{0}^{t}(t-s)^{-\frac{1+n(\frac{1}{\varrho}-\frac{1}{q})}{2}}s^{-(\frac{3+\alpha}{2}-\frac{n}{2\varrho})} ds.
\end{eqnarray*}
Since by hypothesis we have $\max\{1, \frac{n}{3+\alpha}\}<\varrho<\frac{n}{1+\alpha}$ and $\frac{1}{q}-\frac{1}{n}<\frac{1}{\varrho}<\frac{1}{q}$ then, by the same computations as before we have $\displaystyle{\underset{t>0}{\sup}\; t^{\frac{(2+\alpha)q-n}{2q}} \int_{0}^{t}(t-s)^{-\frac{1+n(\frac{1}{\varrho}-\frac{1}{q})}{2}}s^{-(\frac{3+\alpha}{2}-\frac{n}{2\varrho})} ds}<+\infty$, and we can write
\begin{equation}\label{Estimation_ForceExterieureAlpha}
\left\|\int_{0}^{t}\mathfrak{g}_{t-s}\ast f(s,\cdot)ds\right\|_{\mathcal{E}}\leq C\|f\|_{\mathcal{E}_f}.
\end{equation}
\end{itemize}
%%%%%%%%%%%%%%%%%%%%%%%%%%%%%%%%%%%%%%%%%%%%%%%%%%%
With the estimates (\ref{Estimation_DonneeInitialeAlpha}), (\ref{Estimation_BilineaireAlpha}) and (\ref{Estimation_ForceExterieureAlpha}) if the quantity $\|\theta_0\|_{\dot{B}^{-\frac{2q-n}{q},q}_\infty}+\|f\|_{\mathcal{E}_f}$ is small enough, by applying a fixed point argument we obtain a global mild solution of the integral problem (\ref{Formulation_IntegraleAlpha}) and the Theorem \ref{Theoreme3} is now proven. \hfill $\blacksquare$
%%%%%%%%%%%%%%%%%%%%%%%%%%%%%%%%%%%%%%%%%%%%%%%%%%%
%%%%%%%%%%%%%%%%%%%%%%%%%%%%%%%%%%%%%%%%%%%%%%%%%%%
%%%%%%%%%%%%%%%%%%%%%%%%%%%%%%%%%%%%%%%%%%%%%%%%%%%
%%%%%%%%%%%%%%%%%%%%%%%%%%%%%%%%%%%%%%%%%%%%%%%%%%%
\appendix \mysection{Proof of the Proposition \ref{Proposition1}}\label{Secc_ProofProposition1}
In this section we will first establish a pointwise estimate for rough operators of the form (\ref{Def_OperatorAlpha}) and then we will deduce the whished Lebesgue-norm control announced in Proposition \ref{Proposition1}. Recall that these operators are given by the expression 
$$T^\alpha(\phi)(x)=p.v.\int_{\mathbb{R}^n}\frac{\Omega(y/ |y|)}{|y|^{n-\alpha}}\phi(x-y)dy,\qquad 0<\alpha<n-1,$$
where the function $\Omega:\mathbb{S}^{n-1}\longrightarrow \mathbb{R}$ is such that 
$\Omega\in L^1(\mathbb{S}^{n-1})$, $\displaystyle{\int_{\mathbb{S}^{n-1}}\Omega \ d\sigma=0}$ and $\Omega \in L^{\rho}(\mathbb{S}^{n-1})$ with $1<\rho<n$. Similar rough type singular integral operators where studied in \cite{Hoang1}, \cite{Hoang} and \cite{Li} but to the best of our knowledge, the boundedness of this particular operator (\emph{i.e.} where $\Omega\in L^{\rho}(\mathbb{S}^{n-1})$ with $1<\rho<n$ and with $0<\alpha<n-1$) was not studied before. 
%%%%%%%%%%%%%%%%%%%%%%%%%%%%%%%%%%%%%%%%%%%%%%%%%%%
\begin{Remarque}
Note that, although the operator $T^\alpha$ can be defined for all $0<\alpha<n$, we will only establish  here a pointwise bound when $0<\alpha<n-1$. The restriction $0<\alpha<n-1$ seems to be essentially technical but this is enough for our purposes. The case $n-1\leq \alpha<n$ will not be studied here and is left open. 
\end{Remarque}

%%%%%%%%%%%%%%%%%%%%%%%%%%%%%%%%%%%%%%%%%%%%%%%%%%%
In order to perform our computations, we need two ingredients. The first one is the the Hardy-Littlewood maximal operator $\mathscr{M}$ of a (locally integrable) function $\phi:\mathbb{R}^n\longrightarrow \mathbb{R}$ defined by 
$$\mathscr{M}(\phi)(x)=\displaystyle{\underset{B \ni x}{\sup } \;\frac{1}{|B|}\int_{B }|\phi(y)|dy},$$
which is bounded in the Lebesgue spaces $L^p(\mathbb{R}^n)$ with $1<p\leq +\infty$. The second ingredient that we need is the Poincaré-Sobolev inequality: for a regular function $\phi$ and for all ball $B(x,r)$ such that $B(x,r)\subset supp(\phi)$ we have
\begin{equation}\label{PoincareSobolev_inequality}
\left(\frac{1}{|B(x,r)|}\int_{B(x,r)}|\phi(y)-\phi_{B_r}|^{\sigma}dy\right)^{\frac{1}{\sigma}}\leq Cr\left(\frac{1}{|B(x,r)|}\int_{B(x,r)}|\nabla \phi(y)|^{q}dy\right)^{\frac{1}{q}}. 
\end{equation}
for $1\leq q<n$ and $1\leq \sigma\leq \frac{nq}{n-q}$, where $\phi_{B_r}$ stands for the average $\displaystyle{\phi_{B_r}=\frac{1}{B(x,r)}\int_{B(x,r)}\phi(y)dy}$. See a proof of this inequality in \cite[Theorem 3.14]{Kinnunen}.\\

With all these ingredients, we start by defining the operator
\begin{equation}\label{Def_OperatorTstar}
T^{*, \alpha}(\phi)(x)=\underset{t>0}{\sup}\left|\int_{\{|y|>t\}}\frac{\Omega(y/|y|)}{|y|^{n-\alpha}}\phi(x-y)dy\right|,
\end{equation}
and we consider 
$$T^{t, \alpha}(\phi)(x)=\int_{\{|y|>t\}}\frac{\Omega(y/|y|)}{|y|^{n-\alpha}}\phi(x-y)dy,$$
note that we have $T^{*, \alpha}(\phi)(x)=\underset{t>0}{\sup}|T^{t, \alpha}(\phi)(x)|$ and that $|T^{\alpha}(\phi)(x)|\leq T^{*, \alpha}(\phi)(x)$. Now, for some $k_0\in \mathbb{Z}$ so that $2^{k_0-2}<t \leq 2^{k_0-1}$, we write
$$T^{t, \alpha}(\phi)(x)=\int_{\{t<|y|\leq 2^{k_0-1}\}}\frac{\Omega(y/|y|)}{|y|^{n-\alpha}}\phi(x-y)dy+\sum_{k\geq k_0}\int_{\{2^{k-1}<|y|\leq 2^{k}\}}\frac{\Omega(y/|y|)}{|y|^{n-\alpha}}\phi(x-y)dy.$$
Using the fact that the function $\Omega$ is of null integral, we can introduce some constants in the previous expression to obtain
$$T^{t, \alpha}(\phi)(x)=\int_{\{t<|y|\leq 2^{k_0-1}\}}\frac{\Omega(y/|y|)}{|y|^{n-\alpha}}(\phi(x-y)-c_{k_0})dy+\sum_{k\geq k_0}\int_{\{2^{k-1}<|y|\leq 2^{k}\}}\frac{\Omega(y/|y|)}{|y|^{n-\alpha}}(\phi(x-y)-c_k)dy,$$
from which we deduce
\begin{eqnarray*}
|T^{t,\alpha}(\phi)(x)|&\leq &\sum_{k\in \mathbb{Z}}\int_{\{2^{k-1}<|y|\leq 2^{k}\}}\left|\frac{\Omega(y/|y|)}{|y|^{n-\alpha}}(\phi(x-y)-c_k)\right|dy\\
&\leq &C\sum_{k\in \mathbb{Z}}\frac{1}{2^{k(n-\alpha)}}\int_{\{|y|\leq 2^{k}\}}\left|\Omega(y/|y|)(\phi(x-y)-c_k)\right|dy.
\end{eqnarray*}
Now, by the H\"older inequality with $\frac{1}{\rho}+\frac{1}{\rho'}=1$ and $1<\rho<n$, we write
$$|T^{t, \alpha}(\phi)(x)|\leq \sum_{k\in \mathbb{Z}}\frac{2^{k\alpha}}{2^{kn}}\left(\int_{\{|y|\leq 2^{k}\}}|\Omega(y/|y|)|^\rho dy\right)^{\frac{1}{\rho}}\left(\int_{\{|y|\leq 2^{k}\}}|\phi(x-y)-c_k|^{\rho'}dy\right)^{\frac{1}{\rho'}}.$$
Introducing the variable $z=2^{-k} y$, by a change of variables in the first integral above we obtain
$$|T^{t, \alpha}(\phi)(x)|\leq C\sum_{k\in \mathbb{Z}}\frac{2^{k\alpha}}{2^{kn(1-\frac{1}{\rho})}}\left(\int_{\{|z|\leq 1\}}|\Omega(z/|z|)|^\rho dz\right)^{\frac{1}{\rho}}\left(\int_{\{|y|\leq 2^{k}\}}|\phi(x-y)-c_k|^{\rho'}dy\right)^{\frac{1}{\rho'}},$$
and rewriting this formula we have
$$|T^{t, \alpha}(\phi)(x)|\leq C\sum_{k\in \mathbb{Z}}2^{k\alpha}\left(\int_{\{|z|\leq 1\}}|\Omega(z/|z|)|^\rho dz\right)^{\frac{1}{\rho}}\left(\frac{1}{2^{kn}}\int_{\{|y|\leq 2^{k}\}}|\phi(x-y)-c_k|^{\rho'}dy\right)^{\frac{1}{\rho'}}.$$
For the second integral above, we consider now the ball $B(x,2^k)$ and we fix the constant $c_k=\phi_{B_k}=\frac{1}{|B(x,2^k)|}\displaystyle{\int_{B(x,2^k)}\phi(y)dy}$, so we can write (since $\omega_n 2^{kn}=|B(x,2^k)|$, where $\omega_n=|B(0,1)|$ is the volume of the $n$-dimensional unit ball):
$$|T^{t, \alpha}(\phi)(x)|\leq C\sum_{k\in \mathbb{Z}}2^{k\alpha}\left(\int_{\{|z|\leq 1\}}|\Omega(z/|z|)|^\rho dz\right)^{\frac{1}{\rho}}\left(\frac{1}{|B(x,2^k)|}\int_{B(x,2^k)}|\phi(y)-\phi_{B_k}|^{\rho'}dy\right)^{\frac{1}{\rho'}}.$$
We study now more in detail the first integral above, we thus have
\begin{eqnarray*}
|T^{t, \alpha}(\phi)(x)|&\leq &C\left(\int_{0}^1\int_{\mathbb{S}^{n-1}}|\Omega(\xi/|\xi|)|^\rho d\sigma(\xi)r^{n-1}dr\right)^{\frac{1}{\rho}}\sum_{k\in \mathbb{Z}}2^{k\alpha}\left(\frac{1}{|B(x,2^k)|}\int_{B(x,2^k)}|\phi(y)-\phi_{B_k}|^{\rho'}dy\right)^{\frac{1}{\rho'}}\\
&\leq &C\left(\int_{\mathbb{S}^{n-1}}|\Omega(\xi/|\xi|)|^\rho d\sigma(\xi)\right)^{\frac{1}{\rho}}\sum_{k\in \mathbb{Z}}2^{k\alpha}\left(\frac{1}{|B(x,2^k)|}\int_{B(x,2^k)}|\phi(y)-\phi_{B_k}|^{\rho'}dy\right)^{\frac{1}{\rho'}},\end{eqnarray*}
so we obtain
$$|T^{t, \alpha}(\phi)(x)|\leq C \|\Omega\|_{L^\rho(\mathbb{S}^{n-1})}\sum_{k\in \mathbb{Z}}2^{k\alpha}\left(\frac{1}{|B(x,2^k)|}\int_{B(x,2^k)}|\phi(y)-\phi_{B_k}|^{\rho'}dy\right)^{\frac{1}{\rho'}}.$$
Now we apply the Poincaré-Sobolev inequality given in (\ref{PoincareSobolev_inequality}) to obtain
\begin{eqnarray}
|T^{t, \alpha}(\phi)(x)|&\leq &C \|\Omega\|_{L^\rho(\mathbb{S}^{n-1})}\sum_{k\in \mathbb{Z}}2^{k\alpha}\left(\frac{1}{|B(x,2^k)|}\int_{B(x,2^k)}|\phi(y)-\phi_{B_k}|^{\rho'}dy\right)^{\frac{1}{\rho'}}\notag\\
&\leq &C \|\Omega\|_{L^\rho(\mathbb{S}^{n-1})}\sum_{k\in \mathbb{Z}}
2^{k(1+\alpha)}\left(\frac{1}{|B(x,2^k)|}\int_{B(x,2^k)}|\nabla \phi(y)|^{q}dy\right)^{\frac{1}{q}},\label{Somme1}
\end{eqnarray}
where $\frac{n}{n-1}<\rho'$ (since $1<\rho<n$ and $\frac{1}{\rho}+\frac{1}{\rho'}=1$) and where $\rho'\leq \frac{nq}{n-q}$. Note that we thus have $\frac{n}{n-1}<\rho'\leq \frac{nq}{n-q}$ which leads us to the condition $1<\frac{\rho n}{\rho n+\rho-n}\leq q<n$. 
We study now the sum
$$\mathcal{S}=\sum_{k\in \mathbb{Z}}2^{k(1+\alpha)}\left(\frac{1}{|B(x,2^k)|}\int_{B(x,2^k)}|\nabla \phi(y)|^{q}dy\right)^{\frac{1}{q}},$$
and since we can use the following triangle inequality, for all fixed $x\in \mathbb{R}^n$ and $k\in \mathbb{Z}$:
\begin{eqnarray*}
\left(\frac{1}{|B(x,2^k)|}\int_{B(x,2^k)}|\nabla \phi(y)|^{q}dy\right)^{\frac{1}{q}}&\leq &\left(\frac{1}{|B(x,2^k)|}\int_{B(x,2^k)}|\nabla \phi(y)|^{q}\mathds{1}_{\{|x-y|\leq 2^{k-1}\}}dy\right)^{\frac{1}{q}}\\
&&+\left(\frac{1}{|B(x,2^k)|}\int_{B(x,2^k)}|\nabla \phi(y)|^{q}\mathds{1}_{\{2^{k-1}\leq |x-y|\leq 2^k\}}dy\right)^{\frac{1}{q}},
\end{eqnarray*}
we then have
\begin{equation}\label{Somme2}
\begin{split}
\mathcal{S}&\leq \underbrace{\sum_{k\in \mathbb{Z}}2^{k(1+\alpha)}\left(\frac{1}{|B(x,2^k)|}\int_{\{|x-y|\leq 2^{k-1}\}}|\nabla \phi(y)|^{q}dy\right)^{\frac{1}{q}}}_{\mathcal{S}_1} \\
&+\underbrace{\sum_{k\in \mathbb{Z}}2^{k(1+\alpha)}\left(\frac{1}{|B(x,2^k)|}\int_{\{2^{k-1}\leq |x-y|\leq 2^k\}}|\nabla \phi(y)|^{q}dy\right)^{\frac{1}{q}}}_{\mathcal{S}_2}.
\end{split}
\end{equation}
\begin{itemize}
\item For the first term of (\ref{Somme2}) we have
\begin{eqnarray*}
\mathcal{S}_1&=&2^{1+\alpha}\sum_{k\in \mathbb{Z}}2^{(k-1)(1+\alpha)}\left(\frac{1}{|B(x,2^k)|}\int_{\{|x-y|\leq 2^{k-1}\}}|\nabla \phi(y)|^{q}dy\right)^{\frac{1}{q}}\\
&=&2^{(1+\alpha)-\frac{n}{q}}\sum_{k\in \mathbb{Z}}2^{(k-1)(1+\alpha)}\left(\frac{1}{|B(x,2^{k-1})|}\int_{\{|x-y|\leq 2^{k-1}\}}|\nabla \phi(y)|^{q}dy\right)^{\frac{1}{q}},
\end{eqnarray*}
which is 
$$\mathcal{S}_1=2^{(1+\alpha)-\frac{n}{q}}\sum_{k\in \mathbb{Z}}2^{k(1+\alpha)}\left(\frac{1}{|B(x,2^k)|}\int_{B(x,2^k)}|\nabla \phi(y)|^{q}dy\right)^{\frac{1}{q}},$$
and we obtain the formula
\begin{equation}\label{FormulaS1}
\mathcal{S}_1=2^{(1+\alpha)-\frac{n}{q}}\mathcal{S}.
\end{equation}

\item For the second term of (\ref{Somme2}) we have
$$\mathcal{S}_2=\sum_{k\in \mathbb{Z}}2^{k(1+\alpha)}\left(\frac{1}{|B(x,2^k)|}\int_{\{2^{k-1}\leq |x-y|\leq 2^k\}}|\nabla \phi(y)|^{q}dy\right)^{\frac{1}{q}}.$$
Introducing a parameter $0<K<+\infty$ that will be fixed later, we write
\begin{eqnarray}
\mathcal{S}_2&=&\underbrace{\sum_{k\leq  \lfloor \log_2(K)\rfloor }2^{k(1+\alpha)}\left(\frac{1}{|B(x,2^k)|}\int_{\{2^{k-1}\leq |x-y|\leq 2^k\}}|\nabla \phi(y)|^{q}dy\right)^{\frac{1}{q}}}_{(A)}\notag\\
&&+\underbrace{\sum_{k>  \lfloor \log_2(K)\rfloor }2^{k(1+\alpha)}\left(\frac{1}{|B(x,2^k)|}\int_{\{2^{k-1}\leq |x-y|\leq 2^k\}}|\nabla \phi(y)|^{q}dy\right)^{\frac{1}{q}}}_{(B)}.\label{Somme22}
\end{eqnarray}
The term $(A)$ in the expression (\ref{Somme22}) above is treated as follows: 
\begin{eqnarray*}
(A)&=&\sum_{k\leq  \lfloor \log_2(K)\rfloor }2^{k(1+\alpha)}\left(\frac{1}{|B(x,2^k)|}\int_{\{2^{k-1}\leq |x-y|\leq 2^k\}}|\nabla \phi(y)|^{q}dy\right)^{\frac{1}{q}}\\
&\leq &\sum_{k\leq  \lfloor \log_2(K)\rfloor }2^{k(1+\alpha)}\left(\frac{1}{|B(x,2^k)|}\int_{B(x,2^k)}|\nabla \phi(y)|^{q}dy\right)^{\frac{1}{q}},
\end{eqnarray*}
and since we have the control $\displaystyle{\left(\frac{1}{|B(x,2^k)|}\int_{B(x,2^k)} |\nabla \phi(y)|^{q} dy\right)^{\frac{1}{q}}\leq \big(\mathscr{M}(|\nabla \phi|^{q})(x)\big)^{\frac{1}{q}}}$, we obtain
\begin{eqnarray}
(A)&\leq &\sum_{k\leq  \lfloor \log_2(K)\rfloor }2^{k(1+\alpha)}\big(\mathscr{M}(|\nabla \phi|^{q})(x)\big)^{\frac{1}{q}}=\big(\mathscr{M}(|\nabla \phi|^{q})(x)\big)^{\frac{1}{q}} \sum_{k\leq  \lfloor \log_2(K)\rfloor }2^{k(1+\alpha)}\notag \\
&\leq& C K^{1+\alpha}\big(\mathscr{M}(|\nabla \phi|^{q})(x)\big)^{\frac{1}{q}}.\label{SommeA}
\end{eqnarray}
For the term $(B)$ in the expression (\ref{Somme22}) we write
\begin{eqnarray*}
(B)&=&\sum_{k>  \lfloor \log_2(K)\rfloor }2^{k(1+\alpha)}\left(\frac{1}{|B(x,2^k)|}\int_{\{2^{k-1}\leq |x-y|\leq 2^k\})}|\nabla \phi(y)|^{q}dy\right)^{\frac{1}{q}}\\
&\leq &\sum_{k>  \lfloor \log_2(K)\rfloor }2^{k(1+\alpha)}\left(\frac{1}{|B(x,2^k)|}\int_{B(x,2^k)}|\nabla \phi(y)|^{q}dy\right)^{\frac{1}{q}},
\end{eqnarray*}
and we have 
\begin{eqnarray*}
(B)&\leq &C\sum_{k>  \lfloor \log_2(K)\rfloor }2^{k(1+\alpha-\frac{n}{q})}\left(\int_{B(x,2^k)}|\nabla \phi(y)|^{q}dy\right)^{\frac{1}{q}}.
\end{eqnarray*}
Since $\displaystyle{\left(\int_{B(x,2^k)}|\nabla \phi(y)|^{q}dy\right)^{\frac{1}{q}}}\leq \|\nabla \phi\|_{L^q}$, we obtain
\begin{eqnarray*}
(B)\leq C\|\nabla \phi\|_{L^{q}}\sum_{k>  \lfloor \log_2(K)\rfloor }2^{k(1+\alpha-\frac{n}{q})}.
\end{eqnarray*}
But since $1+\alpha-\frac{n}{q}<0$, as we have by hypothesis $1<q<\frac{n}{1+\alpha}$ (recall also that $\alpha<n-1$), we obtain that the previous sum converges and we can write 
\begin{equation}\label{SommeB}
(B)\leq C\|\nabla \phi\|_{L^{q}}K^{1+\alpha-\frac{n}{q}}.
\end{equation}
With the estimates (\ref{SommeA}) and (\ref{SommeB}) at hand, we come back to the expression (\ref{Somme22}) to obtain the inequality
$$\mathcal{S}_2\leq C\left(K^{1+\alpha}\big(\mathscr{M}(|\nabla \phi|^{q})(x)\big)^{\frac{1}{q}}+\|\nabla \phi\|_{L^{q}}K^{1+\alpha-\frac{n}{q}}\right).$$
If we set $K=\left(\frac{\|\nabla \phi\|_{L^{q}}}{\big(\mathscr{M}(|\nabla \phi|^{q})(x)\big)^{\frac{1}{q}}}\right)^{\frac{q}{n}}$, we have
\begin{equation}\label{FormulaS2} 
\mathcal{S}_2\leq C \big(\mathscr{M}(|\nabla \phi|^{q})(x)\big)^{\frac{1}{q}-\frac{(1+\alpha)}{n}}\|\nabla \phi\|_{L^{q}}^{\frac{q(1+\alpha)}{n}}.
\end{equation}
\end{itemize}
With the estimate (\ref{FormulaS1}) for the term $\mathcal{S}_1$ and the inequality (\ref{FormulaS2}) for the term $\mathcal{S}_2$, getting back to the expression (\ref{Somme2}) we have the control:
\begin{eqnarray*}
\mathcal{S}&\leq & \mathcal{S}_1+\mathcal{S}_2\\
&\leq & 2^{(1+\alpha)-\frac{n}{q}}\mathcal{S}+C \big(\mathscr{M}(|\nabla \phi|^{q})(x)\big)^{\frac{1}{q}-\frac{(1+\alpha)}{n}}\|\nabla \phi\|_{L^{q}}^{\frac{q(1+\alpha)}{n}}.
\end{eqnarray*}
Since $q<\frac{n}{1+\alpha}$ we have $2^{(1+\alpha)-\frac{n}{q}}<1$ and we obtain 
$$\mathcal{S}(1- 2^{(1+\alpha)-\frac{n}{q}})\leq C \big(\mathscr{M}(|\nabla \phi|^{q})(x)\big)^{\frac{1}{q}-\frac{(1+\alpha)}{n}}\|\nabla \phi\|_{L^{q}}^{\frac{q(1+\alpha)}{n}},$$
from which we deduce that $\mathcal{S}\leq C\big(\mathscr{M}(|\nabla \phi|^{q})(x)\big)^{\frac{1}{q}-\frac{(1+\alpha)}{n}}\|\nabla \phi\|_{L^{q}}^{\frac{q(1+\alpha)}{n}}$. Thus, coming back to the expression (\ref{Somme1}), we have 
$$|T^{t,\alpha}(\phi)(x)|\leq C \|\Omega\|_{L^\rho(\mathbb{S}^{n-1})}\big(\mathscr{M}(|\nabla \phi|^{q})(x)\big)^{\frac{1}{q}-\frac{(1+\alpha)}{n}}\|\nabla \phi\|_{L^{q}}^{\frac{q(1+\alpha)}{n}},$$
and from this uniform estimate we can obtain the control
$$T^{*, \alpha}(\phi)(x)\leq C \|\Omega\|_{L^\rho(\mathbb{S}^{n-1})}\big(\mathscr{M}(|\nabla \phi|^{q})(x)\big)^{\frac{1}{q}-\frac{(1+\alpha)}{n}}\|\nabla \phi\|_{L^{q}}^{\frac{q(1+\alpha)}{n}},$$
from which we deduce
$$|T^{\alpha}(\phi)(x)|\leq C \|\Omega\|_{L^\rho(\mathbb{S}^{n-1})}\big(\mathscr{M}(|\nabla \phi|^{q})(x)\big)^{\frac{1}{q}-\frac{(1+\alpha)}{n}}\|\nabla \phi\|_{L^{q}}^{\frac{q(1+\alpha)}{n}},$$
 since we have $|T^{\alpha}(\phi)(x)|\leq T^{*, \alpha}(\phi)(x)$. 
 Now, taking the $L^{\mathfrak{s}}$-norm to both sides of this pointwise estimate, we have the inequality
$$\|T^{\alpha}(\phi)\|_{L^{\mathfrak{s}}}\leq C_\Omega  \left\|\big(\mathscr{M}(|\nabla \phi|^{q})\big)^{\frac{1}{q}-\frac{(1+\alpha)}{n}}\right\|_{L^{\mathfrak{s}}}\|\nabla \phi\|_{L^{q}}^{\frac{q(1+\alpha)}{n}},$$
which we rewrite as (using the property $\||f|^\sigma\|_{L^{\mathfrak{s}}}=\|f\|_{L^{\sigma \mathfrak{s}}}^{\sigma}$ for the Lebesgue norms):
$$\|T^{\alpha}(\phi)\|_{L^{\mathfrak{s}}}\leq C_\Omega \left\|\mathscr{M}(|\nabla \phi|^{q})\right\|_{L^{\mathfrak{s}(\frac{1}{q}-\frac{(1+\alpha)}{n})}}^{\frac{1}{q}-\frac{(1+\alpha)}{n}}\|\nabla \phi\|_{L^{q}}^{\frac{q(1+\alpha)}{n}}.$$
Since by hypothesis we have $\mathfrak{s}(\frac{1}{q}-\frac{(1+\alpha)}{n})>1$, thus the maximal function $\mathscr{M}$ is bounded in the Lebesgue space $L^{\mathfrak{s}(\frac{1}{q}-\frac{(1+\alpha)}{n})}(\mathbb{R}^n)$ and we can write
$$\|T^{\alpha}(\phi)\|_{L^\mathfrak{s}}\leq C_\Omega \left\||\nabla \phi|^{q}\right\|_{L^{\mathfrak{s}(\frac{1}{q}-\frac{(1+\alpha)}{n})}}^{\frac{1}{q}-\frac{(1+\alpha)}{n}}\|\nabla \phi\|_{L^{q}}^{\frac{q(1+\alpha)}{n}}\leq C_\Omega \left\|\nabla \phi\right\|_{L^{\mathfrak{s}(1-\frac{q(1+\alpha)}{n})}}^{1-\frac{q(1+\alpha)}{n}}\|\nabla \phi\|_{L^{q}}^{\frac{q(1+\alpha)}{n}},$$
where we used again the property $\||f|^\sigma\|_{L^\mathfrak{s}}=\|f\|_{L^{\sigma \mathfrak{s}}}^{\sigma}$. Since $q=\mathfrak{s}(1-\frac{q(1+\alpha)}{n})$ (recall that $\mathfrak{s}=\frac{qn}{n-q(1+\alpha)}$) we can  obtain the inequality 
$$\|T^{\alpha}(\phi)\|_{L^\mathfrak{s}}\leq C_\Omega \left\|\nabla \phi\right\|_{L^{q}}^{1-\frac{q(1+\alpha)}{n}}\|\nabla \phi\|_{L^{q}}^{\frac{q(1+\alpha)}{n}},$$
which is
$$\|T^{\alpha}(\phi)\|_{L^\mathfrak{s}}\leq C_\Omega \left\|\nabla \phi\right\|_{L^{q}},$$
and this ends the proof of the Proposition \ref{Proposition1}. \hfill$\blacksquare$\\

%%%%%%%%%%%%%%%%%%%%%%%%%%%%%%%%%%%%%%%%%%%%%%%%%%%
%%%%%%%%%%%%%%%%%%%%%%%%%%%%%%%%%%%%%%%%%%%%%%%%%%%
\noindent {\bf Conflict of interest.} We declare that we do not have any commercial or associative interest that represents a conflict of interest in connection with the work submitted.\\

%%%%%%%%%%%%%%%%%%%%%%%%%%%%%%%%%%%%%%%%%%%%%%%%%%%
\noindent {\bf Acknowledgment.} This work was supported by the GDRI ECO-Math.

%%%%%%%%%%%%%%%%%%%%%%%%%%%%%%%%%%%%%%%%%%%%%%%%%%%

%%%%%%%%%%%%%%%%%%%%%%%%%%%%%%%%%%%%%%%%%%%%%%%%%%%%%%%%%%%%%%%%%%%%%%%%%%%%%%%%%%%%%%%%%%%%%%%%%%%%%%%%%%%%%%%%%%%%%%%%%%%%%%%%%%%%%%%%%%%%%%%%%%%%%%%%%%%%%%%%%%%%%%%%%%%%%%%%%%%%%%%%%%%%%%%%%%%%%%%%%%%%
\end{document}